\documentclass[11pt,a4paper
]{article}
\usepackage{amsmath,amssymb,amsthm}
\usepackage[margin=30mm]{geometry}
\usepackage[linktocpage=true]{hyperref}
\usepackage{cleveref}   
\usepackage{tikz}
\usetikzlibrary{arrows.meta,decorations.markings}
\usepackage{pgfplots}

\usepackage[pagewise, mathlines]{lineno} 
\usepackage{time}

\theoremstyle{plain}
\newtheorem{thm}{Theorem}[section]
\newtheorem{prop}[thm]{Proposition}
\newtheorem{lem}[thm]{Lemma}
\newtheorem{cor}[thm]{Corollary}
\theoremstyle{definition}

\newtheorem{rem}[thm]{Remark}
\newtheorem{eg}[thm]{Example}

\theoremstyle{plain}
\newtheorem{thmA}{Theorem}

\numberwithin{equation}{section}
\newcommand{\bb}[1]{\mathbb{#1}}
\newcommand{\mca}[1]{\mathcal{#1}}
\newcommand{\R}{\bb{R}}
\newcommand{\Z}{\bb{Z}}
\newcommand{\Q}{\bb{Q}}
\newcommand{\C}{\bb{C}}

\newcommand{\ol}{\overline}

\newcommand{\surj}{\twoheadrightarrow}
\newcommand{\spmx}[1]{{\small\begin{pmatrix}#1\end{pmatrix}}}
\DeclareMathOperator{\tr}{tr}
\newcommand{\internal}[1]{\par\medskip\noindent
\fbox{\parbox{0.95\linewidth}{\small\textbf{Internal note, not for the
published version.}\\ #1}}\par\medskip}
\renewcommand{\internal}[1]{}  

\DeclareMathOperator{\PSL}{PSL}
\DeclareMathOperator{\SL}{SL}
\DeclareMathOperator{\intr}{int}

\input{mkc1_labels}
\makeatletter
\newcommand{\PI}[1]{\@ifundefined{PI@#1}{\textbf{[I:#1?]}}{\cite[\csname PI@#1\endcsname]{UekiMKC1}}}
\newcommand{\PIr}[2]{\cite[\csname PI@#1\endcsname--\csname PIn@#2\endcsname]{UekiMKC1}}
\newcommand{\PIeq}[1]{\@ifundefined{PI@#1}{\textbf{[I:#1?]}}{\cite[\csname PI@#1\endcsname]{UekiMKC1}}}
\makeatother

\title{\bfseries Modular knots obey the Chebotarev law\\ for compact Lie groups}
\author{\Large Jun Ueki}
\date{}

\begin{document}
\maketitle
\begin{abstract}
The modular knots are the closed orbits of the geodesic flow on the complement of a trefoil in $S^{3}$.
McMullen proved that the closed orbits of a mixing pseudo-Anosov flow on a closed $3$-manifold obey the Chebotarev law for a compact Lie group whose identity component is semisimple.
For the modular knots ordered by the trace, we remove the assumption on the group: 
the cusp, along which the linking numbers grow linearly while the hyperbolic lengths grow logarithmically, 
prevents any character from being a function of the length.
For the generic and for the combinatorial ordering we determine the exceptions: 
by the criterion of Parry and Pollicott, 
the law fails exactly for the characters that are functions of the length, and then 
the Frobenius classes are still equidistributed in the fibers over the torus of these characters and in logarithmic density.
For $\mathrm{SU}(2)$ the traces obey the Sato--Tate law. 
\end{abstract}

{\renewcommand{\thefootnote}{}\footnotetext{\textit{Date:} \today \ \ \now \ (GMT=JST$-9$). }}

{\small \setcounter{tocdepth}{2}\tableofcontents}

\section{Introduction}

This article is a sequel to \cite{UekiMKC1}.
There we proved that the modular knots, together with the missing trefoil $K_{2,3}$,
obey the Chebotarev law in the sense of Mazur and McMullen in $S^{3}$ for the generic orderings and for the trace ordering,
and we determined when the law holds for the combinatorial ordering.
In those statements the target of the homomorphism is a finite group,
except in \PI{thm.main.general} \textup{(2)},
which allows a compact Lie group whose identity component is semisimple for the generic orderings.
Here we replace the finite group by an arbitrary compact Lie group.

The Chebotarev law has a version for a compact Lie group $G$, as in McMullen \cite[Section 1]{McMullen2013CM}:
a homomorphism $\rho:\pi_1(M-L_n)\to G$ is required to have dense image, and the conclusion is the equality
\begin{equation} \label{eq.law.compact}
\lim_{N\to \infty}\frac{1}{N}\sum_{i=n+1}^{N}F\bigl(\rho([K_i])\bigr)
=\int_G F(g)\,dg
\end{equation}
for every continuous class function $F\in C(G,\C)$, where $dg$ denotes the Haar probability measure.
For a finite group $G$, the displayed condition recovers the law of \PI{dfn.Chebotarev}.
The \emph{logarithmic density} version of \eqref{eq.law.compact} is the equality
\begin{equation} \label{eq.law.log}
\lim_{N\to \infty}\frac{1}{\log N}\sum_{i=n+1}^{N}
\frac{F\bigl(\rho([K_i])\bigr)}{i}
=\int_G F(g)\,dg ,
\end{equation}
which follows from \eqref{eq.law.compact} by partial summation, but not conversely.

McMullen proved \eqref{eq.law.compact} for the closed orbits of a topologically mixing pseudo-Anosov flow on a closed $3$-manifold, ordered by length,
under the assumption that the identity component of $G$ is semisimple \cite[Theorem 1.4]{McMullen2013CM},
and he showed that the law can fail for $G=S^{1}$ \cite[Section 6]{McMullen2013CM}.
If the flow is a weak-mixing Anosov flow on a compact manifold $M$ and $\rho$ is a character of $\pi_1(M)$,
then such a failure is conjectured not to occur \cite[p.~15]{KatsudaSunada1990PMIHES}.
For the modular knots, the semisimplicity can be dropped for the trace ordering,
and for the other two orderings the obstruction is determined.

We use the three orderings of the modular knots from \cite{UekiMKC1}, recalled in \Cref{subsec.setting.modular}:
the trace ordering by $|\tr\gamma|$,
the generic orderings by the length $\ell_h$ attached to a non-arithmetic roof function $h$,
and the combinatorial ordering by the length $\ell_{\rm c}$ of the cyclic word.
In the first two, $K_{2,3}$ is inserted at a fixed position, and knots of equal length are ordered arbitrarily.

\begin{thmA} \label{thmA.compact}
Let $G$ be a compact Lie group, let $n\in \Z_{\geq 1}$,
let $\rho\colon \pi_1(S^{3}-L_n)\to G$ be a homomorphism with dense image,
where $L_n$ is the union of the first $n$ knots of the ordering under consideration,
and let $F\in C(G,\C)$ be a continuous class function.
\begin{enumerate}
\item The trace ordering satisfies \eqref{eq.law.compact}.
\item Let $h$ be a non-arithmetic roof function, and let $(K_i)_i$ be the generic ordering by $\ell_h$.
Then \eqref{eq.law.compact} holds for every $F$ if and only if there is no pair of a nontrivial continuous character $\chi\colon G\to S^{1}$ and $a\in \R$ such that
\begin{equation} \label{eq.exceptional}
\chi\bigl(\rho([C_\gamma])\bigr)=e^{ia\,\ell_h(\gamma)}
\end{equation}
holds for all but finitely many modular knots $C_\gamma\not\subset L_n$.
\item For the combinatorial ordering, with $L_n$ a finite union of modular
knots together with $K_{2,3}$ or not, the equality
\[\lim_{N\to\infty}
\frac{1}{\#\{\gamma\mid \ell_{\rm c}(\gamma)\leq N\}}
\sum_{\ell_{\rm c}(\gamma)\leq N}F\bigl(\rho([C_\gamma])\bigr)
=\int_G F(g)\,dg,\]
where $\gamma$ runs through the primitive hyperbolic classes with $C_\gamma\not\subset L_n$,
holds for every $F$ if and only if there is no continuous character $\chi\colon G\to\{\pm1\}$ with $\chi(\rho([C_\gamma]))=(-1)^{\operatorname{lk}(C_\gamma,K_{2,3})}$ for all but finitely many $\gamma$.
\item For the trace ordering and for the generic orderings, the
logarithmic density version \eqref{eq.law.log} holds for every $F$.
For the combinatorial ordering, the averages over the levels converge:
\[\lim_{m_0\to\infty}\frac{1}{m_0}\sum_{m\leq m_0}
\frac{1}{P_m}\sum_{\ell_{\rm c}(\gamma)=m}F\bigl(\rho([C_\gamma])\bigr)
=\int_G F(g)\,dg,\]
where $L_n$ is as in \textup{(3)}, $\gamma$ runs through the classes with $C_\gamma\not\subset L_n$,
and $P_m$ denotes the number of those $\gamma$ with $\ell_{\rm c}(\gamma)=m$.
\end{enumerate}
\end{thmA}

For a finite group $G$, assertion \textup{(1)} is \PI{thm.arith.general}, assertion \textup{(3)} is \PI{thm.combinatorial},
and assertion \textup{(2)} recovers \PI{thm.modular},
since an exception in \textup{(2)} requires a character of infinite image; see the beginning of \Cref{sec.compact.law}.
In \textup{(1)} and \textup{(2)}, $L_n$ may be replaced by any finite union $L$ of modular knots, together with $K_{2,3}$ or not:
such an $L$ is contained in some $L_n$, a homomorphism of $\pi_1(S^{3}-L)$ pulls back to $\pi_1(S^{3}-L_n)$ through the surjection of fundamental groups without changing the Frobenius classes of the knots outside $L_n$,
and the finitely many knots in $L_n-L$ do not affect the limits.

Assertion \textup{(2)} is the case of the modular knots of a general criterion of Parry and Pollicott \cite[Chapter 8]{ParryPollicott1990},
which we state in the setting of the general theorem of \PI{thm.main.general}:
a flow with a Markov section with escape locus,
and its closed orbits ordered by the length attached to a non-arithmetic roof function.
The second assertion below is new.

\begin{thmA} \label{thmA.general}
In the setting of \Cref{subsec.setting.markov}, let $G$ be a compact Lie group,
let $\rho\colon \pi_1(M-L_n)\to G$ be a homomorphism with dense image,
where $L_n$ is the union of the first $n$ closed orbits,
and let $F\in C(G,\C)$ be a continuous class function.
\begin{enumerate}
\item The equality \eqref{eq.law.compact} holds for every $F$ if and only if there is no
pair of a nontrivial continuous character $\chi\colon G\to S^{1}$ and $a\in \R$ such that
\begin{equation} \label{eq.exceptional.general}
\chi\bigl(\rho([K_i])\bigr)=e^{ia\,\ell_h(K_i)}
\end{equation}
holds for all but finitely many $i$.
If such $(\chi,a)$ exists, then $a\neq 0$, and the averages $\frac1N\sum_{i=n+1}^{N}\chi(\rho([K_i]))$ have no limit.
\item The logarithmic density version \eqref{eq.law.log} holds for every
$F$.
\end{enumerate}
\end{thmA}

For a group whose identity component is semisimple the exception cannot occur,
which is how the theorem of McMullen \cite[Theorem 1.4]{McMullen2013CM} follows from \cite[Theorem 8.5]{ParryPollicott1990};
\Cref{thmA.general} \textup{(1)} replaces the assumption on $G$ by the description of the exception.
The exception is the case in which the extension of the flow by $\rho$ has a circle factor,
and it does not affect the logarithmic density.
For the modular knots, the logarithmic density version of the law holds for the trace and for the generic orderings and every compact Lie group,
and for the combinatorial ordering the averages over the levels $\ell_{\rm c}=m$ converge to the Haar integral, although,
in the presence of the parity obstruction, the natural density of each of the two fibers of the parity character oscillates between $\frac13$ and $\frac23$;
see \Cref{thmA.compact} \textup{(4)}.

The exceptions in \Cref{thmA.compact} \textup{(2)} occur.
Let $\Psi(\gamma)$ denote the Rademacher symbol,
which equals $\operatorname{lk}(C_\gamma,K_{2,3})$ by the theorem of Ghys \PI{thm.ghys},
and consider the ordering by $\ell_{\rm c}+\varepsilon\Psi$ of \PI{eg.explicit.ordering},
where $\varepsilon\notin\Q$ and $|\varepsilon|<1$.
The character $\rho_\theta$ of $\pi_1(S^{3}-K_{2,3})$ that sends a meridian of $K_{2,3}$ to $e^{2\pi i\theta}$, with $\theta=(1+\varepsilon)/2$,
has dense image in $S^{1}$ and satisfies $\rho_\theta([C_\gamma])=e^{i\pi(\ell_{\rm c}+\varepsilon\Psi)(\gamma)}$, a function of the length,
so that its averages do not converge; see \Cref{eg.circle.failure}.
The reason why the trace ordering admits no exception is the cusp:
the linking numbers of the classes $R^{q}L$ with the components of $L$ are affine functions of $q$ for large $q$,
while the hyperbolic length grows like $2\log q$,
so that no character can be a function of the hyperbolic length; see \Cref{lem.cusp}.

For $G=\mathrm{SU}(2)$ and the trace ordering,
the values $\frac12\tr\rho([C_\gamma])$ are therefore equidistributed with respect to the measure $\frac{2}{\pi}\sqrt{1-x^{2}}\,dx$ on $[-1,1]$,
which is the measure of the Sato--Tate conjecture for elliptic curves;
see \Cref{cor.satotate} and McMullen \cite[Section 1]{McMullen2013CM}.
Write $\pi_1(S^{3}-K_{2,3})=\langle x,y\mid x^{2}=y^{3}\rangle$.
The homomorphisms $\rho\colon \pi_1(S^{3}-K_{2,3})\to \mathrm{SU}(2)$ with nonabelian image form,
up to conjugation, a single open arc parametrized by $\tr\rho(xy)\in(-\sqrt3,\sqrt3)$,
by Klassen \cite[Theorem 1]{Klassen1991TransAMS}, and all but nine of them have dense image.
Two families lie on the arc: the two classes of homomorphisms that factor through $\pi_1(\Sigma(2,3,7))$,
and the two-dimensional sectors of the Jones representations of the braid group $B_3=\pi_1(S^{3}-K_{2,3})$ at the roots of unity $e^{2\pi i/r}$ with $r\geq 4$,
which have dense image for $r\notin\{4,6,10\}$.
The values $\tr\rho([C_\gamma])$ for the former lie in $\Z[2\cos\frac{2\pi}7]$ and are Galois conjugate to one another; see \Cref{sec.su2}.
For points of the arc with distinct values of $|\tr\rho(xy)|$, the Frobenius classes are independent,
in the sense that their joint distribution along the modular knots is the product of the individual distributions.

\begin{thmA} \label{thmA.joint}
Let $j,m\in \Z_{\geq 1}$, and let $\rho_1,\dots,\rho_j:\pi_1(S^{3}-K_{2,3})\to \mathrm{SU}(2)$ be homomorphisms with dense image for
which the values $|\tr\rho_i(xy)|$, over $1\leq i\leq j$, are pairwise distinct.
Then the values
\[\Bigl(\tfrac12\tr\rho_1([C_\gamma]),\ \dots,\
\tfrac12\tr\rho_j([C_\gamma]),\ \Psi(\gamma)\bmod m\Bigr),\]
taken over the modular knots ordered by $|\tr\gamma|$, and also over the generic ordering by $\ell_h$ for every non-arithmetic roof function $h$,
are equidistributed with respect to the product of $j$ copies of $\frac{2}{\pi}\sqrt{1-x^{2}}\,dx$ and of the uniform measure on $\Z/m\Z$.
\end{thmA}

In particular the two Galois conjugate traces attached to $\Sigma(2,3,7)$ are jointly equidistributed,
and the joint distribution of the Frobenius classes and of the Rademacher symbol modulo $m$ is the product distribution.

The law passes to the preimages under finite branched covers,
in every ordering compatible with the trace ordering of the base,
and in particular to the closed geodesics of the subgroups of finite index of $\PSL_2\Z$;
see \Cref{thm.branched.compact} and \Cref{cor.finiteindex.compact},
whose argument is that of \cite[Theorem B]{UekiSC} with the law for the base supplied by \Cref{thmA.compact} \textup{(1)}.

\subsection*{Organization}
\Cref{sec.setting} recalls the setting and the statements of \cite{UekiMKC1} that we use.
\Cref{sec.criterion} proves \Cref{thmA.general}:
the special characters and the circle factors in \Cref{subsec.special},
the law and its failure in \Cref{subsec.law.failure},
the exceptional case in \Cref{subsec.compact.exceptional},
and the logarithmic density in \Cref{subsec.compact.log}.
\Cref{sec.compact.law} proves \Cref{thmA.compact}: the cusp lemma in \Cref{subsec.cusp.lemma},
then the three orderings in \Crefrange{subsec.compact.arith}{subsec.compact.comb},
and the logarithmic density for the three orderings in \Cref{subsec.compact.log.modular}.
\Cref{sec.su2} describes the homomorphisms into $\mathrm{SU}(2)$,
and \Cref{sec.branched.compact} treats the branched covers.

\section{The setting} \label{sec.setting}

We keep the notation of \cite{UekiMKC1}, and we recall the statements that the proofs below modify.
The statements that are applied without change are cited where they are used.

\subsection{The modular knots and the three orderings}
\label{subsec.setting.modular}

Let $M=T^{1}(\PSL_2\Z\backslash \bb{H}^{2})=S^{3}-K_{2,3}$, where $K_{2,3}$ is the trefoil,
and let $\varphi$ denote the modular flow, the geodesic flow on $M$.
The modular knots are the closed orbits $C_\gamma$ of $\varphi$,
indexed by the primitive hyperbolic conjugacy classes $\gamma$ of $\PSL_2\Z$,
each of which is represented by a positive word in
\[L=\spmx{1&0\\1&1},\qquad R=\spmx{1&1\\0&1},\]
unique up to cyclic permutation, in which both letters occur.
The combinatorial length $\ell_{\rm c}(\gamma)$ is the length of that word,
the hyperbolic length is $\ell(\gamma)=2\log\xi_\gamma$ where $\xi_\gamma+\xi_\gamma^{-1}=|\tr\gamma|$,
and $\ell_h(\gamma)$ is the length attached to a non-arithmetic roof function $h$ as in \Cref{subsec.setting.markov}.
The trace ordering is by $|\tr\gamma|$, equivalently by $\ell$,
the combinatorial ordering is by $\ell_{\rm c}$,
with the knots counted by the threshold $\ell_{\rm c}\leq N$, and a generic ordering is by $\ell_h$.
The Rademacher symbol $\Psi(\gamma)$ is the number of letters $R$ minus the number of letters $L$ in the word of $\gamma$, and by the theorem of Ghys \PI{thm.ghys},
\[\operatorname{lk}(C_\gamma,K_{2,3})=\Psi(\gamma).\]
The Chebotarev law of \PI{dfn.Chebotarev} concerns the sequence of the modular knots in one of the orderings above,
with $K_{2,3}$ inserted at a fixed position, and its version for a compact Lie group is the equality \eqref{eq.law.compact} of the introduction.

We use the following linking numbers with a fixed modular knot \PI{lem.link.cusp}.
Let $\sigma$ be a primitive hyperbolic class,
let $n_R(\sigma)$ and $n_L(\sigma)$ denote the numbers of letters $R$ and $L$ in its word,
and let $m_R(\sigma)$ and $m_L(\sigma)$ denote the largest exponents of $R$ and of $L$ that occur in it.
If $a\geq m_R(\sigma)$ and $R^{a}L\neq\sigma$, then $\operatorname{lk}(C_{R^{a}L},C_\sigma)=n_R(\sigma)$,
and if $b\geq m_L(\sigma)$ and $RL^{b}\neq\sigma$, then $\operatorname{lk}(C_{RL^{b}},C_\sigma)=n_L(\sigma)$.

\subsection{Markov sections with an escape locus}
\label{subsec.setting.markov}

The general theorem \PI{thm.main.general} concerns a flow $\varphi$ on a $3$-manifold $M$ with a Markov section $\mca{R}$ with an escape locus $(\ell^{+},\ell^{-})$ in the sense of \PI{dfn.markov}.
Here $\mca{R}$ is the union of finitely many rectangles $\mca{R}_i=A_i\times B_i$,
not assumed compact, $\ell^{+}$ is a union of finitely many stable arcs,
and $\ell^{-}$ is a union of finitely many unstable arcs.
A point of $\mca{R}$ returns to $\mca{R}$ in positive time if and only if it does not lie on $\ell^{+}$,
and in negative time if and only if it does not lie on $\ell^{-}$,
and the first return map $f$, with return time $r$,
satisfies the Markov conditions on the transition sets $\mca{R}_{ij}=\{x\in \intr\mca{R}_i\mid f(x)\in \intr\mca{R}_j\} =\intr(A_i)\times B_{ij}$.
The section determines a finite directed graph $\mca{G}$, the shift $\Sigma(\mca{G})$ of its two-sided paths,
and the itinerary of each orbit of $\varphi$ that meets $\mca{R}$.
A \emph{roof function} is a positive H\"older function $h$ on $\Sigma(\mca{G})$,
and the length $\ell_h(\tau)$ of a periodic orbit $\tau$ of the shift, and of the closed orbit of $\varphi$ with itinerary $\tau$,
is the sum of $h$ over the points of $\tau$.
The roof function $h$ is \emph{non-arithmetic} if $\ell_h(\tau_1)/\ell_h(\tau_2)\notin\Q$ holds for some pair of periodic orbits that are realized by closed orbits of $\varphi$,
and the ordering by $\ell_h$ is then called a \emph{generic ordering} \PI{subsec.modular.generic}.
The conditions of \PI{subsec.surj} are stated in terms of the following sets.
An arc of $\mca{R}$ is said to be \emph{singular} if it is an edge of some $\mca{R}_i$,
a stable arc $A_i\times\{b\}$ with $b\in B_i$ an endpoint of some $B_{ij}$,
or one of the arcs of $\ell^{+}$ and of $\ell^{-}$.
The \emph{face} of a singular arc $J$ that is not an arc of $\ell^{-}$ is $\bigcup_{x\in J}\varphi_{[0,r(x)]}(x)$,
where $\varphi_{[0,r(x)]}(x)$ is read as $\varphi_{[0,\infty)}(x)$ for $x\in \ell^{+}$,
and the face of an arc $J$ of $\ell^{-}$ is $\varphi_{(-\infty,0]}(J)$.
Let $\mca{W}$ denote the union of the faces, let $N$ denote the union of the orbits that never meet $\mca{R}$,
and let $\mca{U}$ denote the union of the interiors of the rectangles
and of the flow boxes $E_{ij}=\bigcup_{x\in \mca{R}_{ij}}\varphi_{[0,r(x)]}(x)$ of \PI{subsec.boxes}.
The conditions are the following.
\begin{itemize}
\item[\textup{(E0)}] Every compact orbit segment meets $\mca{R}$ in finitely
many points.
\item[\textup{(E1)}] The closure of every face is a properly embedded
surface in $M$.
\item[\textup{(E2)}] The set $N$ is a union of finitely many orbits.
\item[\textup{(E3)}] $M-\mca{U}\subset \mca{W}\cup N$.
\item[\textup{(E4$'$)}] For every unstable arc $J$ of every $\mca{R}_i$,
almost every point of $J$ has a dense forward orbit in $M$,
and for every stable arc $J$ of every $\mca{R}_i$,
almost every point of $J$ has a dense backward orbit in $M$.
\end{itemize}
A periodic point of $\Sigma(\mca{G})$ is realized if some closed orbit of $\varphi$ has it as its itinerary,
and \PI{thm.main.general} assumes that all but finitely many periodic points are realized, each by a unique point of $\mca{R}$,
for the section and for a refinement along every finite union of closed orbits,
that is, for a Markov section with the same escape locus obtained by cutting the rectangles along finitely many stable and unstable arcs, among them the arcs through the points of those orbits,
and that \textup{(E1)} holds for that refinement.
Under these assumptions, and for a non-arithmetic H\"older roof function $h$,
the closed orbits $(K_i)_i$ of $\varphi$ meeting $\mca{R}$, ordered by $\ell_h$, obey the Chebotarev law,
and \PI{thm.main.general} \textup{(2)} proves the law for a compact Lie group whose identity component is semisimple.
For the modular flow the section is constructed from the Farey tessellation \PI{prop.modular},
its refinements are described in \PI{prop.refine},
and exactly two periodic points fail to be realized, by \PI{prop.modular} \textup{(4)}.

Let $\rho\colon \pi_1(M-L_n)\to G$ be a homomorphism to a compact Lie group, where $L_n=\bigcup_{i\leq n}K_i$,
and let $\Sigma(\mca{G}')$ denote the shift of such a refinement along $L_n$.
The flat cocycle $\tilde\rho\colon \Sigma(\mca{G}')\to G$ of \PI{subsec.frobenius} is defined by fixing a path $\theta_v$ in the graph from a base vertex to each vertex $v$,
and by sending a point $x$ to the value of $\rho$ on the image in $M-L_n$ of the loop $\theta_{x_0}\,e_{x_0x_1}\,\theta_{x_1}^{-1}$,
where $e_{x_0x_1}$ is the edge of $x$ and the graph is mapped to $M-L_n$ by the spine of \PI{sec.spine};
the \emph{Frobenius class} $[\tau]$ of a periodic orbit $\tau$ of the shift is the conjugacy class of the product of the values of $\tilde\rho$ along $\tau$, in the order of the itinerary.
If $\tau$ is realized by a closed orbit $K_i$ that avoids the edges of the refined rectangles,
then $[\tau]=\rho([K_i])$ \PI{lem.homotopic},
and all but finitely many closed orbits avoid the edges \PI{lem.orbits.bijection}.
The extension flow $\Sigma(\mca{G}',h,\tilde\rho)$ is the quotient of $\Sigma(\mca{G}')\times\R\times G$ by $(\sigma(x),t,g\tilde\rho(x))\sim(x,t+h(x),g)$,
with the flow $s\cdot[x,t,g]=[x,s+t,g]$ and the commuting left action $g'\cdot[x,t,g]=[x,t,g'g]$ of $G$.
A \emph{circle factor} of a flow is a continuous map $p$ to $S^{1}$ with $p(s\cdot P)=e^{ias}p(P)$ for some $a>0$,
as in \PI{subsec.mixing}; the flow is topologically weak-mixing if it is topologically transitive and has no circle factor.
The theorem of Parry and Pollicott \PI{thm.PP} says that, if the extension flow is weak-mixing,
then the classes $[\tau]$ of the periodic orbits with $\ell_h(\tau)\leq T$ are equidistributed with respect to the Haar measure as $T\to\infty$,
and \PI{lem.ties} passes from that threshold form to the sequence form \eqref{eq.law.compact},
because the periodic orbits of a given length are negligible in number.

\subsection{The transfer operator} \label{subsec.setting.transfer}

For the trace ordering, \PI{sec.arith} replaces the symbolic argument by the transfer operator of Mayer on the induced system of \PI{subsec.induced},
whose states are finitely many intervals $B_u\subset(0,1)$,
indexed by a finite set $F$, with $X=\bigsqcup_{u\in F}B_u$.
A letter $v$ of the induced system starts at one or more states $u$ and arrives at a state $u_v$,
its branch $S_v^{-1}$ maps $B_u$ into $B_{u_v}$,
and the roof $r_v(z)$ is the hyperbolic length of the corresponding orbit segment;
a word $w=v_1\cdots v_n$, read in the order of the forward itinerary,
has the branch $S_w^{-1}=S_{v_1}^{-1}\circ\cdots\circ S_{v_n}^{-1}$, so that the last letter acts first,
and a closed geodesic $\gamma$ coded by a closed path $w$ satisfies $\sum r_{v_i}=\ell(\gamma)$ along the orbit of the fixed point of $S_w^{-1}$.
The cocycle $\alpha(v)\in G$ is the value of $\rho$ on the loop of the letter $v$,
and $\alpha(w)=\alpha(v_1)\cdots\alpha(v_n)$ is the product in the order of the forward itinerary,
which is the Frobenius class of the closed geodesic coded by $w$ \PI{lem.homotopic}.

For a unitary representation $\varrho$ of $G$ on a finite dimensional space and $s\in\C$,
the transfer operator acts on vector valued functions on $X$ by
\[(\mca{L}_{s,\varrho}w)(z)
=\sum_{v}e^{-s\,r_{v}(z)}\,w\bigl(S_{v}^{-1}z\bigr)\,\varrho(\alpha(v)),\]
the sum being over the letters $v$ that start at the state of $z$;
the operator $\mca{L}^{\rm reg}_{s}$ of the left regular representation acts on functions on $X\times G$ by
\[(\mca{L}^{\rm reg}_{s}W)(z,g)
=\sum_{v}e^{-s\,r_{v}(z)}\,W\bigl(S_{v}^{-1}z,\,\alpha(v)g\bigr),\]
and it contains every $\mca{L}_{s,\varrho}$ with $\varrho$ irreducible.

The operators are nuclear on a space of holomorphic functions on a neighborhood of $\ol{X}$ \PI{lem.nuclear},
the $L$-function $L(s,\varrho)=\det(1-\mca{L}_{s+1,\varrho})/\det(1-\mca{L}_{s,\varrho})$ is meromorphic on the half plane $\mathrm{Re}\,s>1/2$ \PI{prop.fredholm},
and \PI{lem.non-vanishing} describes its behavior on the line $\mathrm{Re}\,s=1$ for a finite group $G$.
\Cref{lem.non-vanishing.compact} below treats a compact Lie group.
The passage from the $L$-functions to the counting of closed geodesics is the Tauberian step \PI{lem.hol},
and the law for the trace ordering after the removal of finitely many modular knots is \PI{thm.arith.general}.

\subsection{Finite quotients} \label{subsec.setting.finite}

Two statements about finite quotients are used in the proofs.
The first is the combinatorial input of the combinatorial ordering.
For a finite group $G$, a map $\alpha$ from the edges of $\mca{G}'$ to $G$,
and a unitary representation $\varrho$ of $G$ on $V_\varrho$, let $L_\varrho$ denote the twisted adjacency matrix,
acting on $\bigoplus_{v}V_\varrho$ over the vertices by $(L_\varrho f)(v)=\sum_{(v,w)}\varrho(\alpha(v,w))f(w)$,
so that $\operatorname{Tr}(L_\varrho^{n})=\sum_{\sigma^{n}x=x}\chi_\varrho([x])$.
By \PI{lem.twisted}, if the graph is irreducible, $\varrho$ is irreducible,
the holonomies of the directed loops at a vertex generate $G$,
and $r$ denotes the spectral radius of the adjacency matrix,
then the spectral radius of $L_\varrho$ is at most $r$,
and it equals $r$ if and only if $\varrho$ is a character $\chi$
and there is $\lambda\in S^{1}$ with $\chi([x])=\lambda^{n}$ for every periodic point $x$ of period $n$.
The second is \PI{lem.infinitely.many}: for a link $L$ consisting of $K_{2,3}$ and finitely many modular knots,
a finite group $G$, a surjection $\rho\colon\pi_1(S^{3}-L)\surj G$,
and $g\in G$, there are infinitely many modular knots $C_\gamma\not\subset L$ with $g\in\rho([C_\gamma])$.

\subsection{Compatible orderings of the preimages}
\label{subsec.setting.compatible}

Let $(K_i)_i$ be a sequence of mutually disjoint knots in $S^{3}$ with union $\mca{K}$,
for instance the modular knots and $K_{2,3}$ in one of the three orderings,
and let $\varpi\colon M_\varpi\to S^{3}$ be a finite cover branched along a finite sublink of $\mca{K}$.
Following \cite{UekiSC}, an ordering $(K'_\nu)_{\nu\in\Z_{\geq 1}}$ of the components of $\varpi^{-1}(\mca{K})$ is said to be \emph{compatible} with $(K_i)_i$ if the following two conditions hold,
where a component of $\varpi^{-1}(K_i)$ has \emph{degree} $k$ if it covers $K_i$ with degree $k$.
\begin{itemize}
\item[\textup{(C1)}] For $\nu\in\Z_{\geq 1}$, let $J(\nu)\in\Z_{\geq 0}\cup\{\infty\}$ denote the
supremum of the integers $j\geq 0$ such that every component of degree one over $K_1,\dots,K_j$ is among $K'_1,\dots,K'_\nu$.
The number of the components of degree one among $K'_1,\dots,K'_\nu$ that lie over some $K_i$ with $i>J(\nu)$ is $o(\nu)$.
\item[\textup{(C2)}] The number of the components of degree at least two
among $K'_1,\dots,K'_\nu$ is $o(\nu)$.
\end{itemize}
The \emph{norm ordering} orders the components by the norm $i^{k}$ of a component of degree $k$ over $K_i$, and,
for the trace ordering of the modular knots,
the ordering by the periods of the lifted flow orders them by $k\,\ell(K_i)$;
both are compatible, as shown in \Cref{sec.branched.compact}.

\section{The criterion of Parry and Pollicott} \label{sec.criterion}

Let $G$ be a compact Lie group, and let $dg$ denote the Haar probability measure.
A sequence $(K_i)_i$ of knots in a $3$-manifold $M$ is said to obey the \emph{Chebotarev law for $G$} if,
for every $n\in \Z_{\geq 1}$, every homomorphism $\rho\colon \pi_1(M-L_n)\to G$ with dense image,
and every continuous class function $F\in C(G,\C)$,
the equality \eqref{eq.law.compact} of the introduction holds.
For a finite group $G$ this is the law of \PI{dfn.Chebotarev}.
By the Peter--Weyl theorem the characters of the irreducible unitary representations of $G$ span a dense subspace of the class functions,
so that \eqref{eq.law.compact} holds for every $F$ if it holds for every character,
and for a nontrivial character the right hand side is $0$.

\Cref{thmA.general} is stated in the setting of \PI{thm.main.general}:
a flow $\varphi$ on a $3$-manifold $M$ with a Markov section with escape locus satisfying \textup{(E0)}--\textup{(E4$'$)},
a non-arithmetic roof function $h$, the sequence $(K_i)_i$ of the closed orbits meeting the section, ordered by $\ell_h$,
and the realization and refinement assumptions of that theorem.
For $n\in \Z_{\geq 1}$ put $L_n=\bigcup_{i\leq n}K_i$,
and let $\Sigma(\mca{G}')$ denote the shift of a refinement along $L_n$ given by those assumptions,
with the flat cocycle $\tilde\rho$ attached to a homomorphism $\rho$ of $\pi_1(M-L_n)$ as in \PI{subsec.frobenius},
so that the class $[\tau]$ of a periodic orbit $\tau$ is the product of the values of $\tilde\rho$ along $\tau$.
A periodic orbit realized by a closed orbit $K_i$ that avoids the edges of the refined rectangles has $[\tau]=\rho([K_i])$ by \PI{lem.homotopic};
all but finitely many periodic orbits are realized, by assumption,
and all but finitely many closed orbits avoid the edges, by \PI{lem.orbits.bijection}.
For $T>0$, let $\pi(T)$ denote the number of the $i$ with $\ell_h(K_i)\leq T$.
The base flow of the refined shift is weak-mixing, because the refined graph is irreducible and the flow has no circle factor \PI{lem.refined.nonarithmetic},
so that $\pi(T)\sim e^{\delta T}/(\delta T)$ holds for the topological entropy $\delta>0$, by \cite[Theorem 6.9]{ParryPollicott1990}.
\Cref{thmA.general} is proved in this section.

Assertion \textup{(1)} is due to Adachi and Sunada \cite[Theorems A and D]{AdachiSunada1987JFA},
who characterize the poles of the $L$-functions on the line $\mathrm{Re}(s)=\delta$,
and to Parry and Pollicott \cite[Chapter 8]{ParryPollicott1990}, who prove the equidistribution in the symbolic setting.
In \cite[p.~134]{ParryPollicott1990}, a character $\chi$ of $G$ is called \emph{special} if there is $a\in\R$
such that $\chi([\tau])=e^{ia\ell_h(\tau)}$ holds for every periodic orbit $\tau$ of the shift;
this condition is equivalent to \eqref{eq.exceptional.general}, by \Cref{lem.finitely.many.exceptions} below.
A number $a\in\R$ is called an \emph{eigenfrequency} of a flow on a space $X$
if there is a nonzero continuous function $p\colon X\to\C$ such that $p(s\cdot P)=e^{ias}p(P)$ holds for every $s\in\R$ and every $P\in X$.
Parry and Pollicott construct a homomorphism from the eigenfrequencies of the extension of the flow by $\rho$ to the special characters,
whose kernel consists of the eigenfrequencies of the base \cite[Proposition 8.4]{ParryPollicott1990},
so that for a weak-mixing base, which is the case here,
nontrivial special characters exist if and only if the extension is not weak-mixing.
They also note the oscillation of the averages of a special character \cite[pp.~135--136]{ParryPollicott1990}.
We include the proof in the form that we use,
with the sequence form and with the finitely many exceptions that the escape locus produces, and we add assertion \textup{(2)}.

For a compact Lie group whose identity component is semisimple,
every continuous character has finite image and \eqref{eq.exceptional.general} cannot hold with $a\neq 0$, by \PI{lem.refined.nonarithmetic};
this is how \PI{thm.main.general} \textup{(2)} and the theorem of McMullen \cite[Theorem 1.4]{McMullen2013CM} follow from \cite[Theorem 8.5]{ParryPollicott1990}.
The example of McMullen \cite[Section 6]{McMullen2013CM} for $G=S^{1}$,
a character that equals $e^{2\pi i\ell}$ on a closed orbit of period $\ell$ of a reparametrized suspension flow, is a special character,
and so are the characters of \Cref{eg.circle.failure} below.
Katsuda and Sunada \cite[Proposition 2.1]{KatsudaSunada1990PMIHES} proved that
a weak-mixing Anosov flow on a compact manifold has a nontrivial special character of the first homology group
if and only if the sum of the stable and the unstable bundles is integrable,
and they conjectured that such a character does not exist \cite[p.~15]{KatsudaSunada1990PMIHES}.

\subsection{Special characters and circle factors} \label{subsec.special}

In this subsection we work in the setting of \Cref{thmA.general},
and we write $\Sigma(\mca{G})$ for the refined shift $\Sigma(\mca{G}')$.
For the modular flow, $\Sigma(\mca{G})$ is the symbolic model of \PI{sec.modular} refined along $L_n$ as in \PI{prop.refine}, and the closed orbits are the modular knots.

\begin{lem}[{\PI{lem.finitely.many}}] \label{lem.finitely.many.exceptions}
Let $c\colon \Sigma(\mca{G})\to S^{1}$ be H\"older continuous,
and for a periodic orbit $\tau$ let $c(\tau)$ denote the product of the values of $c$ along $\tau$.
If $c(\tau)=1$ holds for all but finitely many periodic orbits $\tau$, then it holds for every periodic orbit.
\end{lem}

The assumption of \PI{lem.finitely.many} that there is more than one periodic orbit holds here, since the refined graph is irreducible and its shift has infinitely many periodic orbits.

\begin{prop}[{\cite[Proposition 8.4]{ParryPollicott1990}}]
\label{prop.circle.factor}
Let $G$ be a compact Lie group, and let $\rho\colon \pi_1(M-L_n)\to G$ be a homomorphism with dense image.
The flow on $\Sigma(\mca{G},h,\tilde\rho)$ has a circle factor if and only
if there are a nontrivial continuous character $\chi\colon G\to S^{1}$ and $a\in \R$, $a\neq 0$,
such that $\chi([\tau])=e^{ia\ell_h(\tau)}$ holds for every periodic orbit $\tau$ of $\Sigma(\mca{G})$;
equivalently, such that \eqref{eq.exceptional.general} holds for all but finitely many $i$.
\end{prop}

\begin{proof}
We first note that a continuous character $\chi$ of $G$ with $\chi([\tau])=1$ for every periodic orbit $\tau$ is trivial:
the classes $[\tau]$ contain the holonomies of all directed loops at a vertex,
which are dense in $G$ by \PI{lem.dense} \textup{(1)}.

The equivalence of the two conditions on $(\chi,a)$ follows from \Cref{lem.finitely.many.exceptions},
applied to $c=(\chi\circ\tilde\rho)\,e^{-iah}$,
and from the correspondence of the periodic orbits with the closed orbits recalled before \Cref{thmA.general};
a closed orbit meets $\mca{R}$ finitely often and each point of $\mca{R}$ lies in finitely many rectangles,
so a closed orbit realizes only finitely many periodic orbits,
and finitely many closed orbits correspond to finitely many periodic orbits.

Suppose that the flow on $\Sigma(\mca{G},h,\tilde\rho)$ has a circle factor $p$,
so that $p(s\cdot P)=e^{ias}p(P)$ holds for some $a>0$.
As in the proof of \PI{thm.flat}, the function $P\mapsto p(g\cdot P)/p(P)$ is constant along the flow lines for each $g\in G$,
hence constant because the flow is topologically transitive by \PI{lem.dense};
it equals $\chi(g)$ for a continuous character $\chi$ of $G$,
since $G$ acts by the left multiplication $g'\cdot[x,t,g]=[x,t,g'g]$, which commutes with the identification.
Put $q(x)=p([x,0,e])$ for $x\in \Sigma(\mca{G})$, so that $p([x,t,g])=\chi(g)\,e^{iat}\,q(x)$.
The identification $(\sigma(x),0,g\tilde\rho(x))\sim (x,h(x),g)$ yields $\chi(\tilde\rho(x))\,q(\sigma x)=e^{iah(x)}q(x)$,
and iteration along a periodic orbit $\tau$ produces $\chi([\tau])=e^{ia\ell_h(\tau)}$.
If $\chi=1$, then $e^{ia\ell_h(\tau)}=1$ for every $\tau$ with $a\neq 0$,
which contradicts \PI{lem.refined.nonarithmetic}.
Hence $\chi\neq 1$.

Conversely, let $(\chi,a)$ be as in the statement;
replacing $(\chi,a)$ by $(\overline{\chi},-a)$ if necessary, we may assume $a>0$.
The H\"older function $c=(\chi\circ\tilde\rho)\,e^{-iah}$ has $c(\tau)=1$ for every periodic orbit,
so it is a coboundary, $c(x)=q(x)/q(\sigma x)$ for a continuous $q\colon \Sigma(\mca{G})\to S^{1}$,
by the theorem of Liv\v{s}ic for functions with values in $S^{1}$ \cite[Theorem 4]{Livsic1972Izv}.
That theorem is stated for an aperiodic transition matrix.
Let $d$ denote the greatest common divisor of the lengths of the directed loops of $\mca{G}$.
If $d>1$, then the vertices of $\mca{G}$ are divided into $d$ classes that the edges permute cyclically,
the map $\sigma^{d}$ preserves the set $\Sigma_0$ of the sequences $x$ with $x_0$ in one of those classes,
and the transition matrix of $\sigma^{d}$ on $\Sigma_0$ is aperiodic.
We apply the theorem to $\sigma^{d}$ on $\Sigma_0$ and to the function $c(x)c(\sigma x)\cdots c(\sigma^{d-1}x)$,
and we extend $q$ from $\Sigma_0$ to $\Sigma(\mca{G})$ by $q(\sigma x)=q(x)/c(x)$.
Then $p([x,t,g])=\chi(g)\,e^{iat}\,q(x)$ is well defined by the computation above, continuous,
and satisfies $p(s\cdot P)=e^{ias}p(P)$ with $a>0$, so that it is a circle factor.
\end{proof}

\subsection{The law and its failure} \label{subsec.law.failure}

\begin{proof}[Proof of \Cref{thmA.general} \textup{(1)}]
Suppose first that no $(\chi,a)$ as in \eqref{eq.exceptional.general} exists.
By \Cref{prop.circle.factor} the flow on $\Sigma(\mca{G},h,\tilde\rho)$ has no circle factor,
and it is topologically transitive by \PI{lem.dense},
so it is topologically weak-mixing, which is the assumption of \PI{thm.PP}.
That theorem proves \eqref{eq.law.compact} in the threshold form for the periodic orbits of the shift, \PI{lem.ties} passes to the sequence form,
and the finitely many periodic orbits that are not realized,
as well as the finitely many closed orbits removed in the proof of \PI{thm.main.general}, do not affect the limit.

Suppose now that $(\chi,a)$ exists.
Then $a\neq 0$, by \Cref{lem.finitely.many.exceptions} and the first two paragraphs of the proof of \Cref{prop.circle.factor},
and $\chi(\rho([K_i]))=e^{ia\ell_h(K_i)}$ holds for all but finitely many $i$.
We show that \eqref{eq.law.compact} fails for the class function $F=\chi$, whose integral is $0$.
Integration by parts yields, for $T\geq 1$,
\[\sum_{\ell_h(K_i)\leq T}e^{ia\ell_h(K_i)}
=e^{iaT}\pi(T)-ia\int_1^{T}e^{iax}\pi(x)\,dx+O(1)
=\frac{e^{(\delta+ia)T}}{\delta T}\Bigl(1-\frac{ia}{\delta+ia}\Bigr)
+o\Bigl(\frac{e^{\delta T}}{T}\Bigr),\]
because $\int_1^{T}e^{(\delta+ia)x}\frac{dx}{\delta x} =\frac{e^{(\delta+ia)T}}{\delta(\delta+ia)T}+O(e^{\delta T}/T^{2})$
and because the error $\pi(x)-e^{\delta x}/(\delta x)=o(e^{\delta x}/x)$ contributes $o(e^{\delta T}/T)$.
Hence
\begin{equation} \label{eq.oscillation}
\frac{1}{\pi(T)}\sum_{\ell_h(K_i)\leq T}\chi\bigl(\rho([K_i])\bigr)
=\frac{\delta}{\delta+ia}\,e^{iaT}+o(1),
\end{equation}
whose modulus tends to $\delta/\sqrt{\delta^{2}+a^{2}}>0$ while its argument rotates,
so that the averages have no limit; this is the asymptotic noted in \cite[pp.~135--136]{ParryPollicott1990}.
These threshold averages are the averages of \eqref{eq.law.compact} at $N=\pi(T)$, so \eqref{eq.law.compact} fails as well.
\end{proof}

\subsection{The exceptional case} \label{subsec.compact.exceptional}

When special characters exist, the Frobenius classes are still equidistributed,
in the fibers of the quotient of $G$ by the common kernel of the special characters,
over a point of the quotient that is a function of the length.
We work in the setting of \Cref{thmA.general}.
Let $\operatorname{Hom}(G,S^{1})$ denote the group of the continuous characters $\chi\colon G\to S^{1}$, which is finitely generated,
and let $S\subset \operatorname{Hom}(G,S^{1})$ denote the set of the special characters,
that is, of the $\chi$ for which some $a\in \R$ satisfies \eqref{eq.exceptional.general};
it is a subgroup, and the number $a$ is determined by $\chi$:
two values $a,a'$ would imply $e^{i(a-a')\ell_h(\tau)}=1$ for every periodic orbit $\tau$,
by \Cref{lem.finitely.many.exceptions} applied to $c=e^{-i(a-a')h}$,
against \PI{lem.refined.nonarithmetic}.
Write $a(\chi)$ for it; then $\chi\mapsto a(\chi)$ is a homomorphism $S\to \R$, and it is injective,
because $a(\chi)=0$ forces $\chi=1$ by \Cref{lem.finitely.many.exceptions} and the first paragraph of the proof of \Cref{prop.circle.factor}.
Since $S$ is finitely generated, the finitely many exceptions in \eqref{eq.exceptional.general} can be taken independent of $\chi\in S$.
Let $G_S=\bigcap_{\chi\in S}\ker\chi$, a closed normal subgroup with compact abelian quotient $A=G/G_S$,
whose character group is $S$ by Pontryagin duality,
since a subgroup of the discrete group $\operatorname{Hom}(G,S^{1})$ is its own double annihilator.
Let $\eta\colon \R\to A=\operatorname{Hom}(S,S^{1})$ denote the continuous homomorphism $t\mapsto (\chi\mapsto e^{ia(\chi)t})$,
so that $\chi(\eta(t))=e^{ia(\chi)t}$ for every $\chi\in S$;
its image is dense in $A$, because a character $\chi\in S$ that is trivial on $\eta(\R)$ has $a(\chi)=0$, hence $\chi=1$.
For a continuous class function $F$ on $G$, let $\ol{F}\colon A\to \C$ denote the average of $F$ over the fibers of $G\to A$,
with respect to the Haar measure of $G_S$; it is a continuous function with $\int_A\ol{F}=\int_G F$.

\begin{prop} \label{prop.exceptional.fibers}
In the setting of \Cref{thmA.general}, for every continuous class function $F$ on $G$,
\[\lim_{T\to\infty}\frac{1}{\pi(T)}
\sum_{\ell_h(K_i)\leq T}\Bigl(F\bigl(\rho([K_i])\bigr)
-\ol{F}\bigl(\eta(\ell_h(K_i))\bigr)\Bigr)=0 .\]
In words: the Frobenius class of $K_i$ lies in the fiber of $G\to A$ over $\eta(\ell_h(K_i))$, up to finitely many $i$,
and the classes are equidistributed in the fibers.
\end{prop}

\begin{proof}
Let $\varrho$ be an irreducible unitary representation of $G$ with character $\chi_\varrho$, and let $V$ be its space.
If $\varrho$ is trivial on $G_S$, then $\varrho$ factors through the abelian group $A$, so $\varrho$ is a special character $\chi$,
and $\chi(\rho([K_i]))=e^{ia(\chi)\ell_h(K_i)}=\chi(\eta(\ell_h(K_i))) =\ol{\chi}(\eta(\ell_h(K_i)))$ for all but finitely many $i$, so the assertion holds for $F=\chi$ exactly.
If $\varrho$ is not trivial on $G_S$, then the subspace $V^{G_S}$ of the $G_S$-invariant vectors is $G$-invariant, because $G_S$ is normal,
hence $V^{G_S}=0$, and $\int_{G_S}\varrho(gg_1)\,dg_1=\varrho(g)\int_{G_S}\varrho(g_1)\,dg_1=0$ for every $g\in G$, so that $\ol{\chi_\varrho}=0$.
It remains to show $\sum_{\ell_h(K_i)\leq T}\chi_\varrho(\rho([K_i]))=o(\pi(T))$ for such $\varrho$.
By the results of Adachi and Sunada \cite[Theorems 4-6 and 5-5]{AdachiSunada1987JFA}, recalled in \PI{subsec.PP},
the $L$-function $L(s,\varrho)$ extends to a nowhere vanishing meromorphic function on a neighborhood of $\mathrm{Re}(s)\geq \delta$.
If $\dim\varrho\geq 2$, then $L(s,\varrho)$ is holomorphic there.
If $\varrho$ is a character $\chi$, then $L(s,\chi)$ has a pole at $s=\delta+it$ if and only if
$\chi([\tau])=e^{it\ell_h(\tau)}$ holds for every periodic orbit $\tau$,
that is, if and only if $\chi$ is a special character with $a(\chi)=t$.
Since $\varrho$ is not trivial on $G_S$, it is not a special character.
Hence $L(s,\varrho)$ is analytic and nonvanishing on a neighborhood of the line $\mathrm{Re}(s)=\delta$.
The same holds for the dual representation of $\varrho$, which is not trivial on $G_S$ either.
The Tauberian argument of \cite[Theorem 8.5]{ParryPollicott1990},
which uses the weak-mixing of the base through the simple pole of the zeta function at $s=\delta$ and the absence of other singularities on the line,
then proves the required estimate: the function $d\,\vartheta_{\mathbf 1}(T)+\mathrm{Re}(u\,\vartheta_{\chi_\varrho}(T))$ is non-decreasing for a unimodular $u$,
where $d=\dim\varrho$ and $\vartheta$ denotes the sum of $\ell_h$ times the character over the periodic orbits and their multiples,
and the Wiener--Ikehara theorem applied to it and to $\vartheta_{\mathbf 1}$ yields $\vartheta_{\chi_\varrho}(T)=o(e^{\delta T})$,
hence the estimate after the removal of the multiples.
The assertion for a general $F$ follows by the Peter--Weyl theorem,
since $F\mapsto \ol{F}$ is continuous for the uniform norms.
\end{proof}

\begin{rem}[{\cite{Sharp2018MathZ}}] \label{rem.sharp}
A related statement, in which the closed orbits are restricted to a homology class, is due to Sharp.
Consider a transitive Anosov flow on a compact manifold such that every homology class contains a closed orbit,
and a topologically transitive extension of the flow by a compact connected Lie group $G$.
Let $[\gamma]$ denote the holonomy class in $G$ of a closed orbit $\gamma$, and let $\ell(\gamma)$ denote its period.
If $G$ is semisimple,
then the classes $[\gamma]$ of the closed orbits in a fixed homology class are equidistributed with respect to the Haar measure
\cite[Theorem 1.2]{Sharp2018MathZ}.
If $G$ is not semisimple, then the same conclusion holds under the following assumption:
if a continuous character $\chi\colon G\to S^{1}$ and $a\in\R$ satisfy $\chi([\gamma])=e^{ia\ell(\gamma)}$
for every null-homologous closed orbit $\gamma$, then $\chi=1$ and $a=0$ \cite[Theorem 1.3]{Sharp2018MathZ}.
This assumption is the analogue of the absence of a pair $(\chi,a)$ as in \eqref{eq.exceptional.general}, for the null-homologous orbits.
The homology class of a modular knot in $S^{3}-K_{2,3}$ is its linking number with $K_{2,3}$.
\end{rem}

\subsection{Logarithmic density} \label{subsec.compact.log}

The failure of \eqref{eq.law.compact} in \Cref{thmA.general} \textup{(1)} and in \Cref{thmA.compact} \textup{(3)} is an oscillation with mean zero,
which disappears in the logarithmic averages.

\begin{lem} \label{lem.log.partial}
Let $(c_i)_{i\geq 1}$ be a bounded sequence of complex numbers, and put $B(N)=\frac1N\sum_{i\leq N}c_i$.
Then
\[\sum_{i\leq N}\frac{c_i}{i}=\sum_{i\leq N}\frac{B(i)}{i}+O(1).\]
In particular, if $B(N)\to 0$, then $\sum_{i\leq N}c_i/i=o(\log N)$.
\end{lem}

\begin{proof}
Abel summation yields $\sum_{i\leq N}c_i/i=\sum_{i<N}iB(i)\bigl(\frac1i-\frac1{i+1}\bigr)+B(N) =\sum_{i<N}\frac{B(i)}{i+1}+B(N)$,
and $\sum_{i<N}B(i)\bigl(\frac1i-\frac1{i+1}\bigr)=O(1)$.
\end{proof}

\begin{proof}[Proof of \Cref{thmA.general} \textup{(2)}]
It suffices to treat $F=\chi_\varrho$ for a nontrivial irreducible unitary representation $\varrho$ of $G$,
and to show $\sum_{2\leq i\leq N}\chi_\varrho(\rho([K_i]))/i=o(\log N)$.
By the proof of \Cref{prop.exceptional.fibers},
either $\sum_{\ell_h(K_i)\leq T}\chi_\varrho(\rho([K_i]))=o(\pi(T))$,
in which case \Cref{lem.log.partial} and \PI{lem.ties} prove the assertion,
or $\varrho$ is a special character $\chi$, so that $\chi(\rho([K_i]))=e^{iaT_i}$ for all but finitely many $i$, where $T_i=\ell_h(K_i)$ and $a\neq 0$.
The sequence $(T_i)_i$ is non-decreasing, and $\pi(T_i^{-})\leq i\leq \pi(T_i)$,
where $\pi(T)\sim e^{\delta T}/(\delta T)$ and the jump $\pi(T_i)-\pi(T_i^{-})$ is $o(\pi(T_i))$ by \PI{lem.ties}.
Taking logarithms yields $\delta T_i-\log(\delta T_i)=\log i+o(1)$,
hence $\delta T_i=\log i+\log\log i+o(1)$, and
\[e^{iaT_i}=\exp\bigl(ib(\log i+\log\log i)\bigr)\,(1+o(1)),
\qquad b=a/\delta .\]
The function $f(x)=\exp(ib(\log x+\log\log x))/x$ has $f'(x)=O(1/x^{2})$,
so that $\sum_{i\leq N}f(i)=\int_{2}^{N}f(x)\,dx+O(1) =\int_{\log 2}^{\log N}e^{ib\psi(u)}\,du+O(1)$ with $\psi(u)=u+\log u$.
Since $\psi'(u)=1+1/u\geq 1$ and $\psi''(u)=O(1/u^{2})$, integration by parts yields
\[\int_{\log 2}^{\log N}e^{ib\psi(u)}\,du
=\Bigl[\frac{e^{ib\psi(u)}}{ib\psi'(u)}\Bigr]_{\log 2}^{\log N}
+\int_{\log 2}^{\log N}\frac{e^{ib\psi(u)}\psi''(u)}{ib\psi'(u)^{2}}\,du
=O(1).\]
Hence $\sum_{i\leq N}\chi(\rho([K_i]))/i=O(1)+o(\log N)$, as required.
\end{proof}

\begin{rem} \label{rem.log.density}
For the primes, the natural density and the logarithmic density agree in the Chebotarev density theorem,
while they differ in the finer questions of the races between residue classes,
where the logarithmic density is the one for which the limits exist,
by Rubinstein and Sarnak \cite{RubinsteinSarnak1994},
under the generalized Riemann hypothesis and the linear independence over $\Q$ of the imaginary parts of the zeros;
their density is that of a set of real numbers $x$ with respect to $dx/x$,
and since $\log i$ is $\delta\ell_h(K_i)$ up to an error $O(\log\log i)$,
the weight $1/i$ in \eqref{eq.law.log} is of the same type, a uniform average over the length.
For the modular knots and a compact Lie group,
the two densities already differ at the level of the law itself, for the generic orderings,
where the natural density does not exist in the presence of a special character while the logarithmic density is the Haar measure.
\end{rem}

\section{The modular knots} \label{sec.compact.law}

\Cref{thmA.compact} is proved in this section.

No assumption is made on the identity component of $G$.
For a compact Lie group whose identity component is semisimple, every continuous character has finite image,
and \textup{(2)} reduces to \PI{thm.main.general} \textup{(2)};
the new content of \textup{(1)} and \textup{(2)} concerns the case in which the identity component of $G$ is not semisimple,
and an exception in \textup{(2)} requires a continuous character of $G$ of infinite image.
The exceptions depend on the ordering.
For the trace ordering there is no exception.
For a generic ordering exceptions occur, and \Cref{prop.exceptional.explicit} determines them for the orderings of \PI{eg.explicit.ordering};
see also \Cref{eg.circle.failure}.
For the combinatorial ordering the only exception is the parity of the linking number with $K_{2,3}$, as in \PI{thm.combinatorial}.

\subsection{The cusp lemma} \label{subsec.cusp.lemma}

Assertion \textup{(1)} follows from the next lemma:
a character that is a function of the hyperbolic length is trivial,
because the linking numbers of the classes $R^{q}L$ with a fixed knot are affine functions of $q$ for large $q$, while the lengths grow logarithmically.
A parabolic element was used in the same way by Dal'bo and Peign\'e \cite[Lemme VII.10]{DalboPeigne1996AIF},
and by Babillot and Peign\'e \cite[Proposition 2.2]{BabillotPeigne2000ASENS} for the characters of the homology.

\begin{lem} \label{lem.generate.linking}
Let $L$ be a finite union of modular knots, together with $K_{2,3}$ or not,
and let $c$ denote the number of components of $L$.
Then, for every finite set $E$ of modular knots,
the classes in $H_1(S^{3}-L;\Z)\cong \Z^{c}$ of the modular knots $C_\gamma\not\subset L$ that are not in $E$ generate $H_1(S^{3}-L;\Z)$.
\end{lem}

\begin{proof}
Let $H\subset \Z^{c}$ denote the subgroup generated by the classes of the modular knots $C_\gamma\not\subset L$ outside $E$.
If $H\neq \Z^{c}$, then the finitely generated abelian group $\Z^{c}/H$ has a quotient $\Z/p\Z$ for some prime $p$,
so there is a surjection $\rho_p\colon \pi_1(S^{3}-L)\surj \Z/p\Z$,
the composite of the abelianization with a surjection $\Z^{c}\surj \Z/p\Z$ that kills $H$.
By \PI{lem.infinitely.many}, applied to $\rho_p$ pulled back to $\pi_1(S^{3}-(L\cup K_{2,3}))$ through the surjection of fundamental groups if $K_{2,3}\not\subset L$,
every nonzero element of $\Z/p\Z$ is the class $\rho_p([C_\gamma])$ of infinitely many modular knots,
in particular of one outside $E$, which is a contradiction.
\end{proof}

\begin{lem}[the cusp lemma] \label{lem.cusp}
Let $L$ be a finite union of modular knots, together with $K_{2,3}$ or not, let $G$ be a compact Lie group,
let $\rho\colon \pi_1(S^{3}-L)\to G$ be a homomorphism with dense image,
let $\chi\colon G\to S^{1}$ be a continuous character, and let $a\in \R$.
Assume that $\chi(\rho([C_\gamma]))=e^{ia\,\ell(\gamma)}$ holds for all but finitely many modular knots $C_\gamma\not\subset L$, where $\ell$ is the hyperbolic length.
Then $a=0$ and $\chi=1$.
\end{lem}

\begin{proof}
Let $L'$ denote the union of the modular knots in $L$,
and let $C_{\sigma_1},\dots,C_{\sigma_m}$ denote its components.
The character $\chi\circ\rho$ of $\pi_1(S^{3}-L)$ factors through $H_1(S^{3}-L;\Z)$,
which is free abelian on the meridians of the components of $L$,
and the class of a loop $\gamma$ is $\sum_{K}\operatorname{lk}(\gamma,K)\,\mu_K$,
the sum being over the components $K$ of $L$ with meridians $\mu_K$.
Hence there are $\theta_0,\theta_1,\dots,\theta_m\in \R$,
with $\theta_0=0$ if $K_{2,3}\not\subset L$, such that
\[\chi\bigl(\rho([C_\gamma])\bigr)
=\exp 2\pi i\Bigl(\theta_0\operatorname{lk}(C_\gamma,K_{2,3})
+\sum_{j=1}^{m}\theta_j\operatorname{lk}(C_\gamma,C_{\sigma_j})\Bigr)\]
holds for every modular knot $C_\gamma\not\subset L$.

Put $\gamma_q=R^{q}L$ for $q\in \Z_{\geq 1}$.
By \PI{thm.ghys}, $\operatorname{lk}(C_{\gamma_q},K_{2,3})=\Psi(\gamma_q)=q-1$,
and by \PI{lem.link.cusp}, $\operatorname{lk}(C_{\gamma_q},C_{\sigma_j})=n_R(\sigma_j)$ holds
for every $q\geq \max_j m_R(\sigma_j)$ with $\gamma_q\notin\{\sigma_1,\dots,\sigma_m\}$.
The matrix $R^{q}L$ has trace $q+2$, so that $\ell(\gamma_q)=2\operatorname{arccosh}\frac{q+2}{2}$,
and $\ell(\gamma_{q+1})-\ell(\gamma_q)\to 0$ as $q\to\infty$.
For all large $q$ the assumption implies
\[e^{2\pi i\theta_0}
=\frac{\chi(\rho([C_{\gamma_{q+1}}]))}{\chi(\rho([C_{\gamma_q}]))}
=e^{ia(\ell(\gamma_{q+1})-\ell(\gamma_q))}.\]
The left hand side does not depend on $q$ and the right hand side tends to $1$, so both sides equal $1$:
then $a(\ell(\gamma_{q+1})-\ell(\gamma_q))\in 2\pi\Z$ for all large $q$,
which forces $a=0$ because the differences are positive and tend to $0$.
Therefore $\chi(\rho([C_\gamma]))=1$ holds for all but finitely many modular knots $C_\gamma\not\subset L$.
By \Cref{lem.generate.linking}, the classes of those knots generate $H_1(S^{3}-L;\Z)$,
so that $\chi\circ\rho$ is trivial on $\pi_1(S^{3}-L)$.
Since $\rho$ has dense image and $\chi$ is continuous, $\chi=1$.
\end{proof}

\subsection{The trace ordering} \label{subsec.compact.arith}

We prove the following version of \PI{lem.non-vanishing} for a compact group.

\begin{lem} \label{lem.non-vanishing.compact}
In the setting of \Cref{subsec.setting.transfer}, let $G$ be a compact Lie group, and assume the following two conditions.
\begin{itemize}
\item[\textup{(N1$'$)}] For any two states $u$ and $u'$, the values $\alpha(w)$ of the words $w$ whose last letter starts at $u$ and whose first letter arrives at $u'$ are dense in $G$.
\item[\textup{(N2$'$)}] If $\chi\colon G\to S^{1}$ is a continuous character
and $a\in \R$ satisfies $\chi(\alpha(w))=e^{ia\ell(\gamma)}$
for every closed geodesic $\gamma$ coded by a closed path $w$ of the system, then $a=0$ and $\chi=1$.
\end{itemize}
Let $s$ satisfy $\mathrm{Re}\,s=1$, let $\varrho$ be an irreducible unitary representation of $G$,
and suppose that $1$ is an eigenvalue of $\mca{L}_{s,\varrho}$.
Then $s=1$, the representation $\varrho$ is trivial,
and the eigenfunction is a constant multiple of a positive fixed point of $\mca{L}_{1,\mathbf 1}$, which is unique up to a positive scalar.
\end{lem}

The case $\chi=1$ of condition \textup{(N2$'$)} follows from condition \textup{(N2)} of \PI{lem.non-vanishing}:
it says that no $a\neq 0$ has $a\ell(\gamma)\in 2\pi\Z$ for all $\gamma$,
which holds if two lengths have an irrational ratio.

\begin{proof}
The proof follows that of \PI{lem.non-vanishing}.
Let $s=1+it$, and let $w$ be a nonzero solution of $\mca{L}_{s,\varrho}w=w$.
For a column vector $\xi$, the function $W(z,g)=w(z)\varrho(g)\xi$ is continuous on $\ol{X}\times G$ and satisfies $\mca{L}^{\rm reg}_{s}W=W$.
We fix $\xi$ such that $W$ is not identically zero.

Since $|e^{-itr_{v}}|=1$, the triangle inequality yields
\begin{equation} \label{eq.domination.compact}
|W|=\bigl|\mca{L}^{\rm reg}_{s}W\bigr|\leq \mca{L}^{\rm reg}_{1}|W| \qquad \text{on } X\times G.
\end{equation}
Let $\mu$ denote the product of Lebesgue measure on $X$ with the Haar measure of $G$.
As in the proof of \PI{lem.non-vanishing}, the change of variables formula yields $\int_{B_u}e^{-r_v(z)}V(S_v^{-1}z)\,dz=\int_{S_v^{-1}(B_u)}V(y)\,dy$ for every state $u$ and every letter $v$ starting at $u$,
and the images $S_v^{-1}(B_u)$ tile $X$ when the pair $(u,v)$ runs over the pairs with $v$ starting at $u$.
The substitution $g\mapsto\alpha(v)g$ preserves the Haar measure.
Hence $\int\mca{L}^{\rm reg}_{1}V\,d\mu=\int V\,d\mu$ holds for every $V\in L^{1}(\mu)$,
so that equality holds in \eqref{eq.domination.compact} at $\mu$-almost every point, hence everywhere by continuity.

We use the following consequence of \textup{(N1$'$)}:
for every point $(z_0,g_0)$ of $X\times G$, the points $(S_{w}^{-1}z_0,\alpha(w)g_0)$,
where $w$ runs through the words whose last letter starts at the state of $z_0$, are dense in $X\times G$.
Indeed, let $[c]$ be a cylinder, that is, the set of the points of $X$ whose forward itinerary begins with a given word $c$, and let $U\subset G$ be a nonempty open set.
A word $w=cw'$ that begins with $c$ satisfies $\alpha(w)=\alpha(c)\,\alpha(w')$, and the point $S_{w}^{-1}z_0$ lies in $[c]$.
By \textup{(N1$'$)}, the word $w'$ can be chosen so that its last letter starts at the state of $z_0$ and $\alpha(c)\,\alpha(w')\,g_0\in U$.
Since every nonempty open subset of $X$ contains a cylinder, the density follows.

The function $|W|$ vanishes nowhere.
Indeed, suppose that $|W|(z_0,g_0)=0$.
By the equality in \eqref{eq.domination.compact} and the positivity of the weights, the function $|W|$ vanishes at $(S_{v}^{-1}z_0,\alpha(v)g_0)$ for every letter $v$ that starts at the state of $z_0$,
and hence, by iteration, at $(S_{w}^{-1}z_0,\alpha(w)g_0)$ for every word $w$ whose last letter starts at the state of $z_0$.
These points are dense, so that $W\equiv 0$ holds by continuity, which is a contradiction.
The same argument applies to every non-negative solution of $V=\mca{L}^{\rm reg}_{1}V$ that is continuous on $\ol{X}\times G$,
because the equation holds on the closed cells as well, by continuity, and the same iteration applies to a zero at a boundary point of a cell.

Write $W=|W|\,\omega$ with $|\omega|=1$.
Equality in \eqref{eq.domination.compact} means that, for the letters $v$ starting at the state of $z$,
the values $W(S_{v}^{-1}z,\alpha(v)g)$ multiplied by the weights $e^{-s\,r_{v}(z)}$ are non-negative multiples of a common unimodular number.
After division by $|W|(z,g)$ this becomes
\begin{equation} \label{eq.cohom.compact}
e^{-it\,r_{v}(z)}\,\omega\bigl(S_{v}^{-1}z,\alpha(v)g\bigr)=\omega(z,g)
\end{equation}
for every $z\in X$, every $g\in G$, and every letter $v$ starting at the state of $z$.

For $g'\in G$, the function $(z,g)\mapsto\omega(z,gg')/\omega(z,g)$ is continuous on $\ol{X}\times G$.
By \eqref{eq.cohom.compact}, applied to the numerator and to the denominator,
it is invariant under the substitutions $(z,g)\mapsto(S_{v}^{-1}z,\alpha(v)g)$,
because the left multiplication by $\alpha(v)$ commutes with the right multiplication by $g'$.
By the density above and by continuity, this function is constant; let $\chi(g')$ denote its value.
The map $\chi\colon G\to S^{1}$ is a continuous homomorphism.
Putting $g=e$ yields $\omega(z,g)=\chi(g)\,u_0(z)$ with $u_0(z)=\omega(z,e)$, and \eqref{eq.cohom.compact} becomes
\begin{equation} \label{eq.cohom.u}
e^{-it\,r_{v}(z)}\,\chi(\alpha(v))\,u_0(S_{v}^{-1}z)=u_0(z).
\end{equation}
Let $\gamma$ be a closed geodesic coded by a closed path $w$, and let $z_w$ denote the fixed point of $S_{w}^{-1}$.
Iterating \eqref{eq.cohom.u} along $w$ at the point $z_w$ yields $\chi(\alpha(w))=e^{it\ell(\gamma)}$,
because the sum of the roofs along the orbit of $z_w$ is $\ell(\gamma)$.
By \textup{(N2$'$)}, we obtain $t=0$ and $\chi=1$.

Hence $s=1$, and \eqref{eq.cohom.u} says that $u_0(S_{v}^{-1}z)=u_0(z)$ holds for every letter $v$ starting at the state of $z$,
so that $\omega=u_0$ is constant by the density above.
Therefore $|W|$ is a positive continuous fixed point of $\mca{L}^{\rm reg}_{1}$.
Such a fixed point is unique up to a positive scalar.
Indeed, let $V_1$ and $V_2$ be two of them, and put $b=\min_{\ol{X}\times G}V_1/V_2$.
Then $V_1-bV_2$ is a non-negative continuous fixed point of $\mca{L}^{\rm reg}_{1}$ that vanishes at a minimizer,
so that it vanishes identically, and $V_1=bV_2$.
For $g'\in G$, the translate $(z,g)\mapsto|W|(z,gg')$ is again a positive continuous fixed point,
since the right translations commute with $\mca{L}^{\rm reg}_{1}$,
and hence it equals $c(g')\,|W|$ with $c(g')>0$.
The map $c\colon G\to\R_{>0}$ is a continuous homomorphism from a compact group, hence trivial,
and $|W|$ is a function of $z$ alone.

Therefore $w(z)\varrho(g)\xi$ does not depend on $g$.
This holds for every column vector $\xi$, since the argument applies to every $\xi$ for which $W$ is not identically zero.
Hence $w(z)\varrho(g)=w(z)$ holds for every $z$ and every $g\in G$.
Since $w\neq 0$ and $\varrho$ is irreducible, the representation $\varrho$ is trivial.
Then $w$ is a constant multiple of the positive fixed point $|W|$ of $\mca{L}_{1,\mathbf 1}$,
which is unique up to a positive scalar by the argument above.
\end{proof}

\begin{proof}[Proof of \Cref{thmA.compact} \textup{(1)}]
We first treat a homomorphism $\rho$ defined on $\pi_1(S^{3}-K_{2,3})$,
and we remove further modular knots at the end.
Each of \PI{lem.nuclear}, \PI{prop.fredholm},
and \PI{lem.hol} concerns one finite dimensional unitary representation at a time,
and their proofs apply verbatim, with \Cref{lem.non-vanishing.compact} in place of \PI{lem.non-vanishing}.
Condition \textup{(N1$'$)} holds because the values of the cocycle on the directed loops at a vertex are dense in $G$,
by \PI{lem.dense} \textup{(1)} and the convention of \PI{subsec.arith.removing} on the cocycle,
and because concatenation with a fixed path joining two given states translates that dense set by a fixed element of $G$.
Condition \textup{(N2$'$)} is the cusp lemma, \Cref{lem.cusp}:
the closed paths of the system code all but finitely many modular knots $C_\gamma\not\subset L_n$,
and $\alpha(w)$ is the class $\rho([C_\gamma])$ by \PI{lem.homotopic}.

Let $\chi$ be a nontrivial irreducible character, of dimension $d=\chi(e)$, and put
\[\vartheta_{\chi}(T)=\sum_{\substack{\gamma,\ k\geq 1\\
k\ell(\gamma)\leq T}}\ell(\gamma)\,\chi\bigl(\rho([C_\gamma])^{k}\bigr).\]
For a unimodular number $u\in \C$, the function $d\,\vartheta_{\mathbf 1}(T)+\mathrm{Re}\bigl(u\,\vartheta_{\chi}(T)\bigr)$ is non-decreasing,
because $d+\mathrm{Re}\bigl(u\,\chi(g)\bigr)\geq 0$ holds for every $g\in G$, and its Laplace--Stieltjes transform equals
\[-d\,\frac{L'}{L}(s,\mathbf 1)
-\frac12\Bigl(u\,\frac{L'}{L}(s,\chi)
+\ol{u}\,\frac{L'}{L}(s,\ol{\chi})\Bigr),\]
which is holomorphic on $\{\mathrm{Re}\,s\geq 1\}$ except for a simple pole at $s=1$ with residue $d$, by \PI{lem.hol}.
The Wiener--Ikehara theorem, in the form for non-decreasing functions \cite[Appendix I, (I.8)]{ParryPollicott1990},
yields $d\,\vartheta_{\mathbf 1}(T) +\mathrm{Re}\bigl(u\,\vartheta_{\chi}(T)\bigr)\sim d\,e^{T}$ and,
for the trivial character alone, $\vartheta_{\mathbf 1}(T)\sim e^{T}$,
so the two choices $u=1$ and $u=i$ imply $\vartheta_{\chi}(T)=o(e^{T})$.
The argument of \PI{thm.arith}, which removes the terms with $k\geq 2$ and passes to the counting function,
applies because $|\chi|\leq d$ holds on $G$, and it yields
\[\sum_{\ell(\gamma)\leq T}\chi\bigl(\rho([C_\gamma])\bigr)=o\bigl(e^{T}/T\bigr),
\qquad
\#\{\gamma\mid \ell(\gamma)\leq T\}\sim e^{T}/T.\]
Since the characters span a dense subspace of the class functions,
the equality \eqref{eq.law.compact} follows for the ordering by $\ell$,
and the passage to the ordering by $|\tr\gamma|$ is the same as in \PI{thm.arith}.
The removal of $L_n$ is the same as in \PI{thm.arith.general},
with condition \textup{(N1$'$)} for the refined system supplied by \PI{lem.surj} and by \PI{lem.dense} \textup{(1)}, as before.
\end{proof}

The closed geodesics of a subgroup of finite index are treated in \Cref{cor.finiteindex.compact},
by passing to the branched cover as in \PI{thm.finiteindex}.
If no modular knot has been removed, then assertion \textup{(1)} could also be derived from the Selberg trace formula, as noted after \PI{thmA.arith};
for $G=S^{1}$ it follows from the estimate of Sarnak and Mozzochi recalled in \PI{rem.burrin.vonessen}.

\begin{cor} \label{cor.satotate}
Let $\rho:\pi_1(S^{3}-K_{2,3})\to \mathrm{SU}(2)$ be a homomorphism with dense image.
Then the values $\frac12\tr\rho([C_\gamma])$, taken over the modular knots ordered by $|\tr\gamma|$,
are equidistributed with respect to the measure $\frac{2}{\pi}\sqrt{1-x^{2}}\,dx$ on $[-1,1]$.
\end{cor}

\begin{proof}
The continuous class functions of $\mathrm{SU}(2)$ are the continuous functions of $\frac12\tr$,
and the push-forward of the Haar measure under $\frac12\tr$ is $\frac{2}{\pi}\sqrt{1-x^{2}}\,dx$, by the Weyl integration formula.
\end{proof}

\Cref{fig.numerics} (a) shows the distribution of the values $\frac12\tr\rho([C_\gamma])$ over the modular knots with $\tr\gamma\leq 6000$ for the homomorphism $\rho$ with $\tr\rho(xy)=\frac12$ of \Cref{rem.numerics.su2}.

\subsection{The generic orderings} \label{subsec.compact.generic}

We apply \Cref{sec.criterion} to the modular flow,
and we determine the special characters for the orderings of \PI{eg.explicit.ordering}.

\begin{proof}[Proof of \Cref{thmA.compact} \textup{(2)}]
This is \Cref{thmA.general} \textup{(1)} for the modular flow on $M=S^{3}-K_{2,3}$,
in which the closed orbits meeting the section are the modular knots and \eqref{eq.exceptional.general} is \eqref{eq.exceptional}, with two adjustments.
If $K_{2,3}\not\subset L_n$, then $\rho$ is pulled back to $\pi_1(S^{3}-(L_n\cup K_{2,3}))$ through the surjection of fundamental groups,
which does not change the Frobenius classes of the modular knots;
and the insertion of $K_{2,3}$ at a fixed position shifts the indices by one, which does not affect the limits.
\end{proof}

The exceptions can be written down when $h$ depends on the first letter only, as in \PI{eg.explicit.ordering}, and no modular knot has been removed.

\begin{prop} \label{prop.exceptional.explicit}
Let $L_n=K_{2,3}$, let $h$ take the value $\alpha>0$ on the sequences with first letter $L$
and $\beta>0$ on those with first letter $R$, with $\alpha/\beta\notin \Q$,
so that $\ell_h=\alpha n_L+\beta n_R$,
and let $\rho\colon \pi_1(S^{3}-K_{2,3})\to G$ have dense image.
Every continuous character $\chi$ of $G$ satisfies $\chi\circ\rho=e^{2\pi i\theta\Psi}$ for a unique $\theta\in \R/\Z$,
and a nontrivial $\chi$ satisfies \eqref{eq.exceptional} for some $a\in \R$ if and only if
\[\theta\in \Z+\Z\,\frac{\alpha}{\alpha+\beta},\qquad \theta\notin \Z,\]
in which case $a=2\pi m/(\alpha+\beta)$ for the unique $m\in \Z\smallsetminus\{0\}$ with $\theta\equiv -m\alpha/(\alpha+\beta)\pmod 1$.
Consequently, if $\Theta\subset \R/\Z$ denotes the subgroup of the $\theta$ with $e^{2\pi i\theta\Psi}=\chi\circ\rho$ for some continuous character $\chi$ of $G$,
then \eqref{eq.law.compact} fails for the ordering by $\ell_h$ if and only if $\Theta$ contains a nonzero element of $(\Z+\Z\,\alpha/(\alpha+\beta))/\Z$.
For $G=S^{1}$ and $\rho=e^{2\pi i\theta\Psi}$ with $\theta\notin \Q$,
the law fails if and only if $\theta\in \Q+\Q\,\alpha/(\alpha+\beta)$.
\end{prop}

\begin{proof}
Since $H_1(S^{3}-K_{2,3};\Z)\cong \Z$ is generated by the meridian and $\operatorname{lk}(C_\gamma,K_{2,3})=\Psi(\gamma)$ by \PI{thm.ghys},
every character of $\pi_1(S^{3}-K_{2,3})$ is of the form $e^{2\pi i\theta\Psi}$, and $\Psi=n_R-n_L$ on the cyclic word.
Condition \eqref{eq.exceptional} means that $a(\alpha n_L+\beta n_R)-2\pi\theta(n_R-n_L)\in 2\pi\Z$ for all pairs $(n_L,n_R)$ of positive integers with $n_L+n_R$ large,
since every such pair is realized by the primitive word $L^{n_L}R^{n_R}$.
Taking differences in $n_L$ and in $n_R$ yields $a\alpha+2\pi\theta\in 2\pi\Z$ and $a\beta-2\pi\theta\in 2\pi\Z$,
and their sum is $a(\alpha+\beta)=2\pi m$ with $m\in \Z$,
and then $\theta\equiv -m\alpha/(\alpha+\beta)\pmod 1$;
conversely these two conditions imply \eqref{eq.exceptional}.
If $m=0$, then $a=0$ and $\theta\in \Z$, so $\chi\circ\rho=1$ and $\chi=1$ by density;
hence $m\neq 0$ for a nontrivial $\chi$, and $\theta\notin \Z$.
This proves the first assertion.
The second assertion is the first one combined with \Cref{thmA.compact} \textup{(2)}:
the law fails if and only if some character satisfies the first assertion,
and a character is trivial if and only if its $\theta$ is $0$, because $\rho$ has dense image.
For $G=S^{1}$ and $\rho=e^{2\pi i\theta\Psi}$ we have $\Theta=\Z\theta$,
and $k\theta\in \Z+\Z\,\alpha/(\alpha+\beta)$ with $k\neq 0$ and $k\theta\notin \Z$ holds for some $k$ if and only
if $\theta\in \Q+\Q\,\alpha/(\alpha+\beta)$ and $\theta\notin \Q$, because $\alpha/(\alpha+\beta)\notin \Q$.
\end{proof}

\begin{eg} \label{eg.circle.failure}
Let $G=S^{1}$, and let $\rho_\theta$ be the character that sends a meridian of $K_{2,3}$ to $e^{2\pi i\theta}$,
so that $\rho_\theta([C_\gamma])=e^{2\pi i\theta\Psi(\gamma)}$;
it has dense image if and only if $\theta\notin \Q$.
Consider the ordering by $\ell_h=\ell_{\rm c}+\varepsilon\Psi$ of \PI{eg.explicit.ordering},
that is, $\alpha=1-\varepsilon$ and $\beta=1+\varepsilon$ with $\varepsilon\notin \Q$, $|\varepsilon|<1$,
and put $\theta=(1+\varepsilon)/2$, so that $\theta\equiv -\alpha/(\alpha+\beta) \pmod 1$.
Then $\theta\notin \Q$, so $\rho_\theta$ has dense image, but
\[\rho_\theta([C_\gamma])=e^{i\pi(1+\varepsilon)(n_R-n_L)}
=e^{i\pi(n_L+n_R)}\,e^{i\pi\varepsilon(n_R-n_L)}
=e^{i\pi\ell_h(\gamma)}\]
holds for every $\gamma$, since $n_R-n_L\equiv n_L+n_R\pmod 2$:
the Frobenius element is a function of the length,
with $a=\pi$ in \eqref{eq.exceptional}, and $m=1$ in \Cref{prop.exceptional.explicit}.
By \eqref{eq.oscillation}, the average of $\rho_\theta([C_\gamma])$ over $\ell_h(\gamma)\leq T$ has modulus tending to $\delta/\sqrt{\delta^{2}+\pi^{2}}$,
where $\delta$ is the exponential growth rate of the number of the classes with $\ell_h\leq T$,
the positive solution of $e^{-\delta\alpha}+e^{-\delta\beta}=1$.
\end{eg}

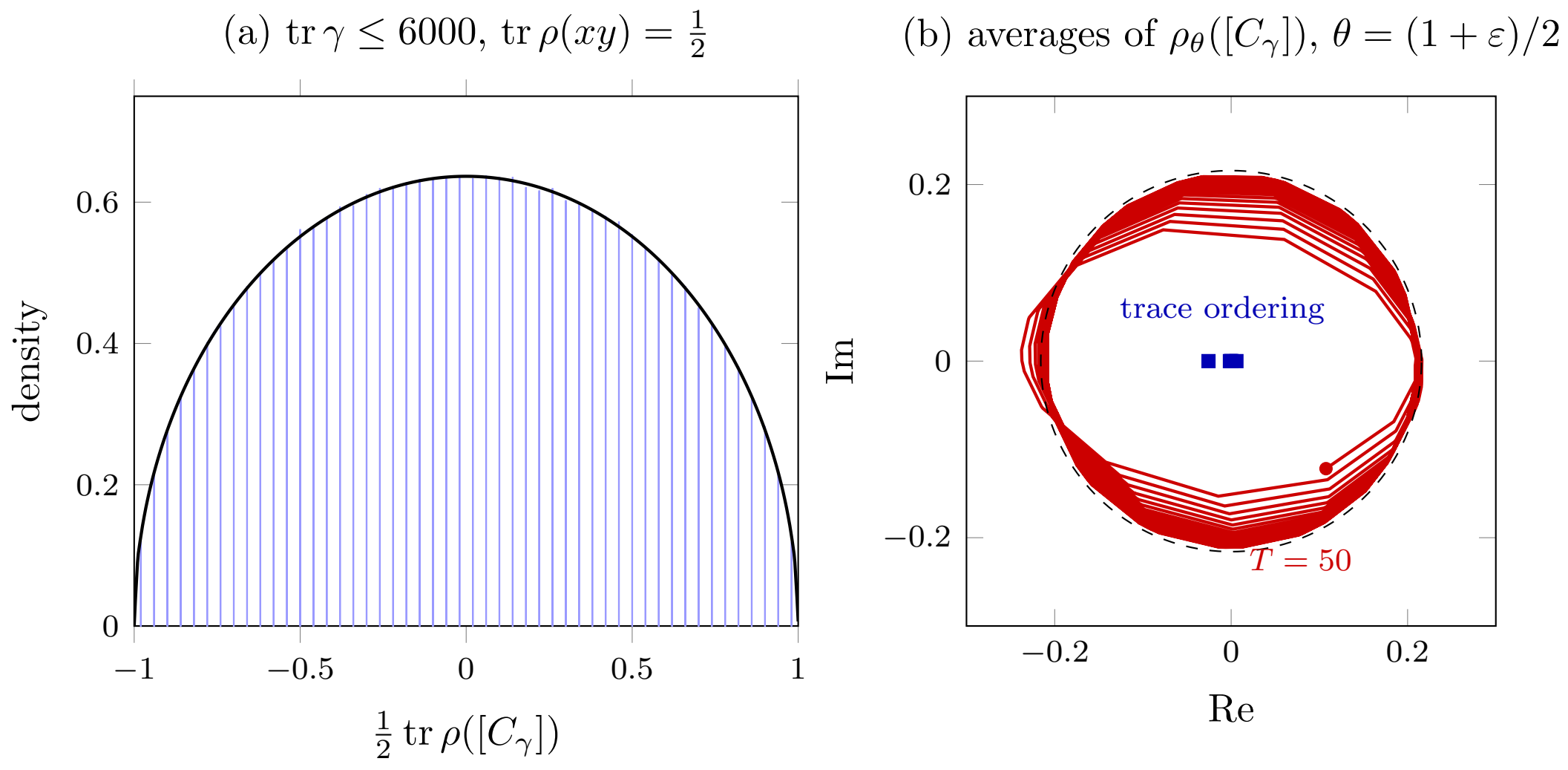
\begin{figure}[t]
\centering
\begin{tikzpicture}
\begin{axis}[name=hist,width=0.5\textwidth,height=0.42\textwidth,
  xmin=-1,xmax=1,ymin=0,ymax=0.75,
  xlabel={$\frac12\tr\rho([C_\gamma])$},ylabel={density},ylabel near ticks,xlabel near ticks,
  title={\small (a) $\tr\gamma\leq 6000$, $\tr\rho(xy)=\frac12$},
  tick label style={font=\scriptsize},label style={font=\small},
  ybar,bar width=0.04,bar shift=0pt]
\addplot[fill=blue!25,draw=blue!40] coordinates {(-0.980,0.1188) (-0.940,0.2133) (-0.900,0.2781) (-0.860,0.3248) (-0.820,0.3613) (-0.780,0.4010) (-0.740,0.4239) (-0.700,0.4550) (-0.660,0.4804) (-0.620,0.4978) (-0.580,0.5160) (-0.540,0.5325) (-0.500,0.5609) (-0.460,0.5633) (-0.420,0.5762) (-0.380,0.5929) (-0.340,0.5946) (-0.300,0.6105) (-0.260,0.6190) (-0.220,0.6220) (-0.180,0.6273) (-0.140,0.6319) (-0.100,0.6337) (-0.060,0.6339) (-0.020,0.6360) (0.020,0.6329) (0.060,0.6369) (0.100,0.6346) (0.140,0.6348) (0.180,0.6204) (0.220,0.6158) (0.260,0.6190) (0.300,0.6019) (0.340,0.5978) (0.380,0.5867) (0.420,0.5796) (0.460,0.5717) (0.500,0.5536) (0.540,0.5352) (0.580,0.5197) (0.620,0.4971) (0.660,0.4803) (0.700,0.4501) (0.740,0.4266) (0.780,0.3985) (0.820,0.3665) (0.860,0.3203) (0.900,0.2819) (0.940,0.2147) (0.980,0.1184)};
\addplot[black,thick,domain=-1:1,samples=200,smooth] {2/pi*sqrt(1-x^2)};
\end{axis}
\begin{axis}[at={(hist.east)},anchor=west,xshift=1.5cm,
  width=0.42\textwidth,height=0.42\textwidth,axis equal image,
  xmin=-0.3,xmax=0.3,ymin=-0.3,ymax=0.3,
  xlabel={$\mathrm{Re}$},ylabel={$\mathrm{Im}$},xlabel near ticks,ylabel near ticks,
  title={\small (b) averages of $\rho_\theta([C_\gamma])$, $\theta=(1+\varepsilon)/2$},
  tick label style={font=\scriptsize},label style={font=\small},
  xtick={-0.2,0,0.2},ytick={-0.2,0,0.2}]
\addplot[red!80!black,thick] coordinates {(0.1075,-0.1218) (0.1844,-0.0685) (0.2088,-0.0222) (0.2098,-0.0180) (0.2102,-0.0003) (0.2053,0.0215) (0.2035,0.0257) (0.1639,0.0790) (0.0602,0.1378) (-0.0770,0.1487) (-0.1821,0.1056) (-0.2287,0.0497) (-0.2296,0.0476) (-0.2376,0.0126) (-0.2372,0.0005) (-0.2345,-0.0115) (-0.2150,-0.0514) (-0.2135,-0.0534) (-0.1423,-0.1144) (-0.0145,-0.1530) (0.1092,-0.1343) (0.1861,-0.0796) (0.1870,-0.0786) (0.2115,-0.0312) (0.2132,-0.0235) (0.2138,0.0006) (0.2072,0.0292) (0.2045,0.0353) (0.1614,0.0925) (0.1605,0.0933) (0.0595,0.1493) (-0.0693,0.1582) (-0.1697,0.1163) (-0.1703,0.1158) (-0.2175,0.0595) (-0.2193,0.0556) (-0.2284,0.0169) (-0.2280,-0.0019) (-0.2246,-0.0175) (-0.2027,-0.0630) (-0.2002,-0.0664) (-0.1269,-0.1290) (-0.1263,-0.1294) (-0.0073,-0.1641) (0.1115,-0.1448) (0.1119,-0.1446) (0.1851,-0.0914) (0.1869,-0.0892) (0.2128,-0.0385) (0.2153,-0.0270) (0.2159,0.0005) (0.2079,0.0351) (0.2038,0.0443) (0.1585,0.1039) (0.1569,0.1054) (0.0587,0.1588) (0.0584,0.1589) (-0.0639,0.1661) (-0.0642,0.1661) (-0.1605,0.1252) (-0.1618,0.1242) (-0.2091,0.0680) (-0.2122,0.0612) (-0.2222,0.0198) (-0.2218,-0.0037) (-0.2179,-0.0225) (-0.1937,-0.0732) (-0.1901,-0.0781) (-0.1177,-0.1397) (-0.1166,-0.1403) (-0.0020,-0.1729) (-0.0018,-0.1729) (-0.0017,-0.1729) (0.1103,-0.1536) (0.1112,-0.1531) (0.1831,-0.1002) (0.1860,-0.0965) (0.2122,-0.0450) (0.2156,-0.0293) (0.2162,0.0016) (0.2070,0.0411) (0.2020,0.0524) (0.1552,0.1136) (0.1525,0.1160) (0.0568,0.1672) (0.0562,0.1674) (0.0561,0.1674) (-0.0605,0.1734) (-0.0605,0.1734) (-0.0611,0.1732) (-0.1539,0.1333) (-0.1562,0.1314) (-0.2034,0.0753) (-0.2078,0.0658) (-0.2184,0.0225) (-0.2180,-0.0061) (-0.2137,-0.0273) (-0.1881,-0.0810) (-0.1835,-0.0874) (-0.1114,-0.1485) (-0.1096,-0.1494) (-0.1096,-0.1495) (-0.0001,-0.1798) (0.0003,-0.1798) (0.0007,-0.1798) (0.1080,-0.1605) (0.1081,-0.1605) (0.1098,-0.1597) (0.1794,-0.1079) (0.1837,-0.1024) (0.2104,-0.0499) (0.2146,-0.0304) (0.2153,0.0025) (0.2052,0.0461) (0.1993,0.0595) (0.1516,0.1215) (0.1481,0.1247) (0.1480,0.1247) (0.0550,0.1739) (0.0539,0.1742) (0.0537,0.1743) (-0.0579,0.1792) (-0.0581,0.1792) (-0.0593,0.1790) (-0.1488,0.1399) (-0.1488,0.1399) (-0.1523,0.1371) (-0.1990,0.0814) (-0.2048,0.0690) (-0.2158,0.0244) (-0.2155,-0.0082) (-0.2108,-0.0315) (-0.1840,-0.0879) (-0.1784,-0.0956) (-0.1784,-0.0956) (-0.1073,-0.1554) (-0.1048,-0.1568) (-0.1047,-0.1569) (0.0009,-0.1855) (0.0015,-0.1855) (0.0023,-0.1855) (0.1051,-0.1664) (0.1052,-0.1664) (0.1079,-0.1650) (0.1759,-0.1142) (0.1759,-0.1142) (0.1817,-0.1068) (0.2085,-0.0541) (0.2134,-0.0311) (0.2142,0.0037) (0.2033,0.0506) (0.1967,0.0657) (0.1967,0.0657) (0.1486,0.1280) (0.1441,0.1320) (0.1441,0.1320) (0.0536,0.1793) (0.0520,0.1798) (0.0515,0.1799) (-0.0556,0.1841) (-0.0560,0.1841) (-0.0579,0.1836) (-0.1444,0.1455) (-0.1445,0.1455) (-0.1494,0.1415) (-0.1955,0.0865) (-0.1955,0.0865) (-0.2026,0.0713) (-0.2139,0.0259) (-0.2136,-0.0105) (-0.2086,-0.0354) (-0.2086,-0.0354) (-0.1810,-0.0935) (-0.1745,-0.1024) (-0.1745,-0.1024) (-0.1043,-0.1611) (-0.1011,-0.1629) (-0.1008,-0.1630) (0.0009,-0.1901) (0.0018,-0.1901) (0.0031,-0.1901) (0.1022,-0.1712) (0.1024,-0.1712) (0.1063,-0.1693) (0.1725,-0.1195) (0.1725,-0.1195) (0.1798,-0.1101) (0.2066,-0.0574) (0.2066,-0.0574) (0.2123,-0.0311) (0.2131,0.0048) (0.2131,0.0048) (0.2017,0.0545) (0.1944,0.0710) (0.1944,0.0710) (0.1462,0.1333) (0.1408,0.1381) (0.1407,0.1382) (0.0526,0.1838) (0.0505,0.1844) (0.0497,0.1846) (-0.0535,0.1882) (-0.0541,0.1881) (-0.0569,0.1874) (-0.1407,0.1503) (-0.1408,0.1502) (-0.1472,0.1450) (-0.1926,0.0908) (-0.1926,0.0907) (-0.2009,0.0729) (-0.2124,0.0270) (-0.2124,0.0270) (-0.2124,0.0270) (-0.2121,-0.0125) (-0.2069,-0.0387) (-0.2069,-0.0387) (-0.1788,-0.0981) (-0.1715,-0.1081) (-0.1714,-0.1081) (-0.1024,-0.1656) (-0.0984,-0.1678) (-0.0979,-0.1681) (0.0005,-0.1938) (0.0018,-0.1939) (0.0038,-0.1939) (0.0994,-0.1753) (0.0996,-0.1751) (0.1048,-0.1726) (0.1694,-0.1239) (0.1694,-0.1238) (0.1782,-0.1125) (0.2049,-0.0602) (0.2049,-0.0602) (0.2049,-0.0602) (0.2113,-0.0308) (0.2122,0.0060) (0.2122,0.0060) (0.2122,0.0060) (0.2003,0.0578) (0.1925,0.0755) (0.1925,0.0756) (0.1443,0.1377) (0.1381,0.1431) (0.1378,0.1433) (0.0521,0.1874) (0.0493,0.1882) (0.0481,0.1885) (-0.0516,0.1916) (-0.0523,0.1914) (-0.0562,0.1905) (-0.1374,0.1542) (-0.1376,0.1541) (-0.1455,0.1476) (-0.1902,0.0943) (-0.1902,0.0943) (-0.1902,0.0942) (-0.1997,0.0738) (-0.2113,0.0277) (-0.2113,0.0277) (-0.2113,0.0277) (-0.2110,-0.0144) (-0.2056,-0.0417) (-0.2056,-0.0417) (-0.2056,-0.0417) (-0.1771,-0.1019) (-0.1691,-0.1128) (-0.1689,-0.1130) (-0.1010,-0.1693) (-0.0963,-0.1719) (-0.0955,-0.1723) (-0.0002,-0.1968) (0.0015,-0.1969) (0.0043,-0.1969) (0.0967,-0.1786) (0.0971,-0.1785) (0.1037,-0.1752) (0.1667,-0.1276) (0.1667,-0.1276) (0.1667,-0.1275) (0.1770,-0.1143) (0.2035,-0.0624) (0.2035,-0.0624) (0.2035,-0.0623) (0.2104,-0.0303) (0.2114,0.0071) (0.2114,0.0071) (0.2114,0.0071) (0.1991,0.0607) (0.1909,0.0794) (0.1908,0.0795) (0.1908,0.0795) (0.1429,0.1412) (0.1359,0.1473) (0.1355,0.1476) (0.0519,0.1903) (0.0485,0.1913) (0.0467,0.1917) (-0.0497,0.1944) (-0.0507,0.1942) (-0.0558,0.1930) (-0.1345,0.1576) (-0.1345,0.1576) (-0.1348,0.1574) (-0.1443,0.1496) (-0.1881,0.0972) (-0.1881,0.0972) (-0.1881,0.0972) (-0.1988,0.0744) (-0.2104,0.0283) (-0.2104,0.0283) (-0.2104,0.0283) (-0.2102,-0.0161) (-0.2046,-0.0442) (-0.2046,-0.0443) (-0.2046,-0.0443) (-0.1758,-0.1050) (-0.1672,-0.1168) (-0.1669,-0.1171) (-0.1669,-0.1171) (-0.1002,-0.1723) (-0.0947,-0.1752) (-0.0936,-0.1757) (-0.0011,-0.1993) (0.0011,-0.1994) (0.0047,-0.1994) (0.0943,-0.1815) (0.0943,-0.1815) (0.0949,-0.1812) (0.1029,-0.1772) (0.1643,-0.1306) (0.1643,-0.1306) (0.1644,-0.1305) (0.1761,-0.1155) (0.2023,-0.0641) (0.2023,-0.0641) (0.2023,-0.0641) (0.2098,-0.0297) (0.2108,0.0081) (0.2108,0.0081) (0.2108,0.0082) (0.1982,0.0631) (0.1896,0.0827) (0.1895,0.0829) (0.1895,0.0829) (0.1418,0.1441) (0.1342,0.1508) (0.1335,0.1513) (0.1335,0.1513) (0.0520,0.1927) (0.0480,0.1938) (0.0456,0.1943) (-0.0478,0.1967) (-0.0478,0.1967) (-0.0492,0.1965) (-0.0555,0.1949) (-0.1320,0.1604) (-0.1320,0.1604) (-0.1324,0.1602) (-0.1434,0.1511) (-0.1864,0.0997) (-0.1864,0.0997) (-0.1865,0.0996) (-0.1981,0.0747) (-0.2098,0.0287) (-0.2098,0.0286) (-0.2098,0.0286) (-0.2096,-0.0177) (-0.2039,-0.0464) (-0.2039,-0.0466) (-0.2039,-0.0466) (-0.1748,-0.1076) (-0.1657,-0.1201) (-0.1653,-0.1205) (-0.1653,-0.1205) (-0.0997,-0.1746) (-0.0935,-0.1779) (-0.0920,-0.1787) (-0.0920,-0.1787) (-0.0021,-0.2014) (-0.0021,-0.2014) (0.0006,-0.2015) (0.0052,-0.2014) (0.0920,-0.1838) (0.0920,-0.1838) (0.0928,-0.1835) (0.1023,-0.1788) (0.1622,-0.1333) (0.1622,-0.1333) (0.1624,-0.1331) (0.1754,-0.1163) (0.2013,-0.0656) (0.2013,-0.0656) (0.2013,-0.0655) (0.2094,-0.0290) (0.2104,0.0090) (0.2104,0.0091) (0.2104,0.0091) (0.1975,0.0652) (0.1886,0.0855) (0.1885,0.0858) (0.1885,0.0858) (0.1411,0.1464) (0.1328,0.1536) (0.1319,0.1543) (0.1319,0.1543) (0.0522,0.1946) (0.0522,0.1946) (0.0477,0.1959) (0.0445,0.1965) (0.0445,0.1965) (-0.0461,0.1986) (-0.0461,0.1986) (-0.0478,0.1984) (-0.0554,0.1965) (-0.1298,0.1628) (-0.1298,0.1628) (-0.1303,0.1625) (-0.1428,0.1522) (-0.1850,0.1017) (-0.1850,0.1017) (-0.1851,0.1016) (-0.1977,0.0747) (-0.2093,0.0289) (-0.2093,0.0289) (-0.2093,0.0288) (-0.2091,-0.0191) (-0.2033,-0.0483) (-0.2033,-0.0485) (-0.2033,-0.0485) (-0.1742,-0.1098) (-0.1645,-0.1229) (-0.1640,-0.1235) (-0.1640,-0.1235) (-0.0994,-0.1765) (-0.0994,-0.1765) (-0.0927,-0.1802) (-0.0906,-0.1811) (-0.0906,-0.1811) (-0.0031,-0.2030) (-0.0031,-0.2030) (0.0000,-0.2032) (0.0057,-0.2031) (0.0057,-0.2031) (0.0900,-0.1859) (0.0900,-0.1859) (0.0911,-0.1854) (0.1020,-0.1800) (0.1605,-0.1355) (0.1605,-0.1355) (0.1607,-0.1352) (0.1750,-0.1168) (0.2006,-0.0668) (0.2006,-0.0667) (0.2006,-0.0667) (0.2090,-0.0283) (0.2101,0.0099) (0.2100,0.0100) (0.2100,0.0100) (0.1971,0.0670) (0.1879,0.0878) (0.1876,0.0882) (0.1876,0.0883) (0.1407,0.1483) (0.1407,0.1483) (0.1318,0.1560) (0.1306,0.1569) (0.1306,0.1569) (0.0527,0.1961) (0.0527,0.1961) (0.0476,0.1976) (0.0436,0.1984) (0.0436,0.1984) (-0.0445,0.2003) (-0.0445,0.2003) (-0.0466,0.1999) (-0.0555,0.1978) (-0.0555,0.1978) (-0.1279,0.1648) (-0.1279,0.1648) (-0.1285,0.1644) (-0.1425,0.1530) (-0.1838,0.1034) (-0.1839,0.1034) (-0.1840,0.1032) (-0.1974,0.0745) (-0.2090,0.0291) (-0.2090,0.0290) (-0.2090,0.0290) (-0.2088,-0.0204) (-0.2029,-0.0500) (-0.2029,-0.0503) (-0.2029,-0.0503) (-0.1737,-0.1115) (-0.1737,-0.1115) (-0.1636,-0.1252) (-0.1629,-0.1260) (-0.1629,-0.1260) (-0.0994,-0.1781) (-0.0994,-0.1781) (-0.0920,-0.1820) (-0.0894,-0.1832) (-0.0894,-0.1832) (-0.0042,-0.2044) (-0.0042,-0.2044) (-0.0006,-0.2046) (0.0062,-0.2045) (0.0062,-0.2045) (0.0882,-0.1876) (0.0882,-0.1876) (0.0895,-0.1871) (0.1019,-0.1809) (0.1019,-0.1809) (0.1590,-0.1373) (0.1590,-0.1373) (0.1593,-0.1370) (0.1748,-0.1171) (0.2000,-0.0677) (0.2000,-0.0677) (0.2000,-0.0676) (0.2088,-0.0276) (0.2099,0.0106) (0.2098,0.0108) (0.2098,0.0108) (0.1967,0.0685) (0.1967,0.0685) (0.1873,0.0898) (0.1870,0.0904) (0.1870,0.0904) (0.1405,0.1498) (0.1405,0.1498) (0.1310,0.1580) (0.1294,0.1592) (0.1294,0.1592) (0.0532,0.1974) (0.0532,0.1974) (0.0476,0.1989) (0.0427,0.2000) (0.0427,0.2000) (-0.0430,0.2017) (-0.0430,0.2017) (-0.0454,0.2013) (-0.0556,0.1988) (-0.0556,0.1988) (-0.1262,0.1666) (-0.1262,0.1666) (-0.1270,0.1661) (-0.1423,0.1535) (-0.1423,0.1535) (-0.1829,0.1049) (-0.1829,0.1049) (-0.1830,0.1046) (-0.1973,0.0743) (-0.2087,0.0292) (-0.2087,0.0291) (-0.2087,0.0290) (-0.2086,-0.0215) (-0.2086,-0.0215) (-0.2027,-0.0513) (-0.2025,-0.0518) (-0.2025,-0.0518) (-0.1734,-0.1130) (-0.1734,-0.1130) (-0.1629,-0.1271) (-0.1620,-0.1282) (-0.1620,-0.1282) (-0.0996,-0.1793) (-0.0996,-0.1793) (-0.0916,-0.1835) (-0.0884,-0.1850) (-0.0884,-0.1850) (-0.0052,-0.2055) (-0.0052,-0.2055) (-0.0012,-0.2057) (0.0067,-0.2056) (0.0067,-0.2056) (0.0866,-0.1891) (0.0866,-0.1891) (0.0881,-0.1884) (0.1019,-0.1815) (0.1019,-0.1815) (0.1577,-0.1390) (0.1577,-0.1389) (0.1581,-0.1385) (0.1747,-0.1172) (0.1747,-0.1172) (0.1995,-0.0686) (0.1995,-0.0685) (0.1996,-0.0683) (0.2087,-0.0269) (0.2087,-0.0269) (0.2097,0.0112) (0.2097,0.0115) (0.2097,0.0116) (0.1965,0.0697) (0.1965,0.0697) (0.1869,0.0914) (0.1864,0.0923) (0.1864,0.0923) (0.1404,0.1510) (0.1404,0.1510) (0.1304,0.1596) (0.1284,0.1611) (0.1284,0.1611) (0.0538,0.1984) (0.0538,0.1984) (0.0477,0.2001) (0.0419,0.2013) (0.0419,0.2013) (-0.0416,0.2029) (-0.0417,0.2029) (-0.0444,0.2024) (-0.0559,0.1995) (-0.0559,0.1995) (-0.1248,0.1681) (-0.1248,0.1681) (-0.1257,0.1675) (-0.1423,0.1538) (-0.1423,0.1538) (-0.1821,0.1061) (-0.1821,0.1061) (-0.1823,0.1058) (-0.1972,0.0740) (-0.1972,0.0740) (-0.1972,0.0740) (-0.2086,0.0293) (-0.2086,0.0291) (-0.2086,0.0290) (-0.2084,-0.0224) (-0.2084,-0.0224) (-0.2025,-0.0525) (-0.2023,-0.0531) (-0.2023,-0.0531) (-0.1733,-0.1141) (-0.1733,-0.1141) (-0.1624,-0.1287) (-0.1613,-0.1301) (-0.1613,-0.1301) (-0.0998,-0.1803) (-0.0998,-0.1803) (-0.0914,-0.1848) (-0.0875,-0.1865) (-0.0875,-0.1865) (-0.0063,-0.2065) (-0.0063,-0.2065) (-0.0018,-0.2067) (0.0073,-0.2065) (0.0073,-0.2065) (0.0851,-0.1903) (0.0851,-0.1903) (0.0869,-0.1896) (0.1020,-0.1820) (0.1020,-0.1820) (0.1566,-0.1403) (0.1566,-0.1403) (0.1571,-0.1398) (0.1747,-0.1172) (0.1747,-0.1172) (0.1747,-0.1172) (0.1991,-0.0693) (0.1992,-0.0691) (0.1992,-0.0689) (0.2087,-0.0262) (0.2087,-0.0262) (0.2087,-0.0262) (0.2097,0.0118) (0.2097,0.0122) (0.2097,0.0123) (0.1963,0.0707) (0.1963,0.0707) (0.1866,0.0928) (0.1860,0.0939) (0.1860,0.0939) (0.1404,0.1520) (0.1404,0.1520) (0.1300,0.1609) (0.1276,0.1628) (0.1276,0.1628) (0.0545,0.1992) (0.0545,0.1992) (0.0478,0.2010) (0.0412,0.2024) (0.0412,0.2024) (-0.0404,0.2039) (-0.0404,0.2039) (-0.0435,0.2033) (-0.0563,0.2001) (-0.0563,0.2001) (-0.1235,0.1694) (-0.1235,0.1694) (-0.1245,0.1687) (-0.1424,0.1540) (-0.1424,0.1540) (-0.1424,0.1540) (-0.1814,0.1072) (-0.1815,0.1071) (-0.1817,0.1067) (-0.1973,0.0736) (-0.1973,0.0736) (-0.1973,0.0736) (-0.2085,0.0294) (-0.2085,0.0292) (-0.2086,0.0290) (-0.2084,-0.0233) (-0.2084,-0.0233) (-0.2084,-0.0233) (-0.2024,-0.0535) (-0.2022,-0.0543) (-0.2022,-0.0543) (-0.1732,-0.1151) (-0.1732,-0.1151) (-0.1621,-0.1300) (-0.1607,-0.1317) (-0.1606,-0.1317) (-0.1002,-0.1810) (-0.1002,-0.1810) (-0.0912,-0.1858) (-0.0866,-0.1879) (-0.0866,-0.1879) (-0.0073,-0.2073) (-0.0073,-0.2073) (-0.0023,-0.2075) (0.0078,-0.2073) (0.0078,-0.2073) (0.0838,-0.1914) (0.0838,-0.1914) (0.0858,-0.1906) (0.1023,-0.1823) (0.1023,-0.1823) (0.1023,-0.1823) (0.1556,-0.1415) (0.1557,-0.1415) (0.1562,-0.1409) (0.1748,-0.1170) (0.1748,-0.1170) (0.1748,-0.1170) (0.1989,-0.0699) (0.1989,-0.0697) (0.1990,-0.0694) (0.2087,-0.0256) (0.2087,-0.0256) (0.2087,-0.0256) (0.2097,0.0122) (0.2097,0.0128) (0.2097,0.0129) (0.1963,0.0716) (0.1963,0.0716) (0.1963,0.0716) (0.1864,0.0939) (0.1857,0.0953) (0.1857,0.0953) (0.1405,0.1527) (0.1405,0.1527) (0.1298,0.1620) (0.1268,0.1642) (0.1268,0.1642) (0.0552,0.1999) (0.0552,0.1999) (0.0480,0.2018) (0.0404,0.2034) (0.0404,0.2034) (-0.0392,0.2047) (-0.0392,0.2047) (-0.0426,0.2041) (-0.0567,0.2006) (-0.0567,0.2006) (-0.0567,0.2006) (-0.1223,0.1705) (-0.1224,0.1705) (-0.1236,0.1697) (-0.1426,0.1540) (-0.1426,0.1540) (-0.1426,0.1540) (-0.1809,0.1081) (-0.1810,0.1080) (-0.1812,0.1075) (-0.1973,0.0732) (-0.1973,0.0732) (-0.1973,0.0732) (-0.2085,0.0295) (-0.2085,0.0291) (-0.2086,0.0289) (-0.2084,-0.0240) (-0.2084,-0.0240) (-0.2084,-0.0240) (-0.2024,-0.0543) (-0.2021,-0.0553) (-0.2021,-0.0554) (-0.1732,-0.1158) (-0.1732,-0.1158) (-0.1732,-0.1158) (-0.1618,-0.1311) (-0.1601,-0.1332) (-0.1601,-0.1332) (-0.1006,-0.1816) (-0.1006,-0.1816) (-0.0912,-0.1866) (-0.0859,-0.1890) (-0.0859,-0.1890) (-0.0082,-0.2079) (-0.0082,-0.2079) (-0.0029,-0.2082) (0.0084,-0.2079) (0.0084,-0.2079) (0.0084,-0.2079) (0.0826,-0.1924) (0.0826,-0.1924) (0.0849,-0.1914) (0.1026,-0.1825) (0.1026,-0.1825) (0.1026,-0.1825) (0.1548,-0.1426) (0.1549,-0.1425) (0.1556,-0.1418) (0.1750,-0.1168) (0.1750,-0.1168) (0.1750,-0.1168) (0.1987,-0.0704) (0.1988,-0.0702) (0.1989,-0.0698) (0.2088,-0.0251) (0.2088,-0.0251) (0.2088,-0.0251) (0.2098,0.0126) (0.2097,0.0133) (0.2097,0.0135) (0.1963,0.0723) (0.1963,0.0723) (0.1963,0.0723) (0.1863,0.0949) (0.1854,0.0966) (0.1854,0.0966) (0.1407,0.1533) (0.1407,0.1533) (0.1407,0.1533) (0.1296,0.1629) (0.1262,0.1655) (0.1262,0.1655) (0.0559,0.2004) (0.0559,0.2004) (0.0483,0.2025) (0.0397,0.2042) (0.0397,0.2042) (0.0397,0.2042) (-0.0381,0.2054) (-0.0381,0.2054) (-0.0418,0.2048) (-0.0572,0.2010) (-0.0572,0.2010) (-0.0572,0.2010) (-0.1213,0.1716) (-0.1214,0.1715) (-0.1228,0.1705) (-0.1429,0.1540) (-0.1429,0.1540) (-0.1429,0.1540) (-0.1804,0.1089) (-0.1805,0.1088) (-0.1809,0.1082) (-0.1975,0.0728) (-0.1975,0.0728) (-0.1975,0.0728) (-0.2085,0.0296) (-0.2086,0.0291) (-0.2086,0.0289) (-0.2084,-0.0246) (-0.2084,-0.0246) (-0.2084,-0.0246) (-0.2024,-0.0549) (-0.2021,-0.0562) (-0.2020,-0.0563) (-0.1733,-0.1164) (-0.1733,-0.1164) (-0.1733,-0.1164) (-0.1617,-0.1320) (-0.1597,-0.1344) (-0.1596,-0.1345) (-0.1010,-0.1821) (-0.1010,-0.1821) (-0.1010,-0.1821) (-0.0912,-0.1873) (-0.0851,-0.1900) (-0.0851,-0.1900) (-0.0851,-0.1900) (-0.0092,-0.2085) (-0.0092,-0.2085) (-0.0034,-0.2087) (0.0091,-0.2085) (0.0091,-0.2085) (0.0091,-0.2085) (0.0815,-0.1932) (0.0816,-0.1932) (0.0841,-0.1922) (0.1030,-0.1827) (0.1030,-0.1827) (0.1030,-0.1827) (0.1541,-0.1435) (0.1542,-0.1434) (0.1550,-0.1426) (0.1752,-0.1166) (0.1752,-0.1166) (0.1752,-0.1166) (0.1985,-0.0709) (0.1987,-0.0706) (0.1988,-0.0701) (0.2089,-0.0246) (0.2089,-0.0246) (0.2089,-0.0246) (0.2099,0.0129) (0.2098,0.0138) (0.2098,0.0140) (0.1963,0.0729) (0.1963,0.0729) (0.1963,0.0729) (0.1862,0.0956) (0.1852,0.0977) (0.1851,0.0977) (0.1410,0.1538) (0.1410,0.1538) (0.1410,0.1538) (0.1295,0.1636) (0.1256,0.1666) (0.1256,0.1666) (0.1255,0.1666) (0.0566,0.2008) (0.0566,0.2008) (0.0566,0.2008) (0.0485,0.2030) (0.0390,0.2049) (0.0390,0.2049) (0.0390,0.2049) (-0.0370,0.2061) (-0.0371,0.2061) (-0.0412,0.2053) (-0.0577,0.2013) (-0.0577,0.2013) (-0.0577,0.2012) (-0.1204,0.1724) (-0.1205,0.1724) (-0.1221,0.1713) (-0.1432,0.1539) (-0.1432,0.1539) (-0.1432,0.1539) (-0.1801,0.1097) (-0.1802,0.1095) (-0.1806,0.1088) (-0.1977,0.0725) (-0.1977,0.0725) (-0.1977,0.0725) (-0.2085,0.0297) (-0.2086,0.0291) (-0.2087,0.0288) (-0.2085,-0.0251) (-0.2085,-0.0251) (-0.2085,-0.0251) (-0.2025,-0.0554) (-0.2020,-0.0571) (-0.2020,-0.0572) (-0.1734,-0.1169) (-0.1734,-0.1169) (-0.1734,-0.1169) (-0.1616,-0.1327) (-0.1592,-0.1356) (-0.1592,-0.1356) (-0.1592,-0.1356) (-0.1015,-0.1824) (-0.1015,-0.1824) (-0.1015,-0.1824) (-0.0913,-0.1878) (-0.0845,-0.1909) (-0.0845,-0.1909) (-0.0845,-0.1909) (-0.0100,-0.2089) (-0.0100,-0.2089) (-0.0100,-0.2089) (-0.0039,-0.2092) (0.0097,-0.2089) (0.0097,-0.2089) (0.0097,-0.2089) (0.0805,-0.1940) (0.0806,-0.1940) (0.0834,-0.1928) (0.1034,-0.1827) (0.1034,-0.1827) (0.1034,-0.1827) (0.1535,-0.1444) (0.1536,-0.1442) (0.1545,-0.1432) (0.1754,-0.1164) (0.1754,-0.1164) (0.1754,-0.1164) (0.1984,-0.0713) (0.1986,-0.0709) (0.1987,-0.0704) (0.2090,-0.0241) (0.2090,-0.0241) (0.2090,-0.0241) (0.2100,0.0131) (0.2099,0.0143) (0.2099,0.0145) (0.1964,0.0733) (0.1964,0.0733) (0.1964,0.0733) (0.1863,0.0963) (0.1850,0.0987) (0.1850,0.0987) (0.1849,0.0987) (0.1413,0.1541) (0.1413,0.1541) (0.1413,0.1541) (0.1295,0.1642) (0.1250,0.1676) (0.1250,0.1676) (0.1250,0.1676) (0.0572,0.2011) (0.0572,0.2011) (0.0572,0.2011) (0.0488,0.2034) (0.0384,0.2056) (0.0384,0.2056) (0.0384,0.2056) (-0.0361,0.2066) (-0.0361,0.2066) (-0.0361,0.2066) (-0.0405,0.2058) (-0.0583,0.2015) (-0.0583,0.2015) (-0.0583,0.2015) (-0.1196,0.1732) (-0.1197,0.1732) (-0.1215,0.1719) (-0.1436,0.1538) (-0.1436,0.1538) (-0.1436,0.1538) (-0.1797,0.1103) (-0.1799,0.1101) (-0.1804,0.1093) (-0.1978,0.0721) (-0.1978,0.0721) (-0.1978,0.0721) (-0.2086,0.0298) (-0.2087,0.0290) (-0.2088,0.0286) (-0.2086,-0.0256) (-0.2086,-0.0256) (-0.2086,-0.0256) (-0.2026,-0.0559) (-0.2020,-0.0578) (-0.2020,-0.0578) (-0.2020,-0.0579) (-0.1736,-0.1172) (-0.1736,-0.1172) (-0.1736,-0.1172) (-0.1616,-0.1333) (-0.1589,-0.1366) (-0.1589,-0.1366) (-0.1588,-0.1366) (-0.1019,-0.1827) (-0.1019,-0.1827) (-0.1019,-0.1827) (-0.0914,-0.1883) (-0.0838,-0.1917) (-0.0838,-0.1917) (-0.0838,-0.1917) (-0.0109,-0.2093) (-0.0109,-0.2093) (-0.0109,-0.2093) (-0.0044,-0.2096) (0.0103,-0.2093) (0.0103,-0.2093) (0.0103,-0.2093) (0.0797,-0.1946) (0.0797,-0.1946) (0.0797,-0.1946) (0.0827,-0.1933) (0.1039,-0.1827) (0.1039,-0.1827) (0.1039,-0.1827) (0.1530,-0.1451) (0.1531,-0.1449) (0.1541,-0.1438) (0.1757,-0.1161) (0.1757,-0.1161) (0.1757,-0.1161) (0.1983,-0.0717) (0.1986,-0.0712) (0.1987,-0.0706) (0.2091,-0.0237) (0.2091,-0.0237) (0.2091,-0.0237) (0.2101,0.0132) (0.2100,0.0147) (0.2100,0.0147) (0.2100,0.0149) (0.1965,0.0737) (0.1965,0.0737) (0.1965,0.0737) (0.1863,0.0968) (0.1848,0.0996) (0.1848,0.0996) (0.1848,0.0997) (0.1416,0.1544) (0.1416,0.1544) (0.1416,0.1544) (0.1296,0.1647) (0.1245,0.1685) (0.1245,0.1685) (0.1245,0.1685) (0.0579,0.2014) (0.0579,0.2014) (0.0579,0.2014) (0.0491,0.2038) (0.0377,0.2061) (0.0377,0.2061) (0.0377,0.2061) (-0.0352,0.2071) (-0.0352,0.2071) (-0.0352,0.2071) (-0.0400,0.2063) (-0.0588,0.2016) (-0.0588,0.2016) (-0.0588,0.2016) (-0.1189,0.1739) (-0.1191,0.1739) (-0.1191,0.1739) (-0.1210,0.1725) (-0.1439,0.1536) (-0.1439,0.1536) (-0.1439,0.1536) (-0.1795,0.1109) (-0.1797,0.1106) (-0.1802,0.1097) (-0.1980,0.0717) (-0.1980,0.0717) (-0.1980,0.0717) (-0.2087,0.0300) (-0.2089,0.0289) (-0.2089,0.0289) (-0.2089,0.0285) (-0.2087,-0.0260) (-0.2087,-0.0260) (-0.2087,-0.0260) (-0.2027,-0.0562) (-0.2021,-0.0585) (-0.2021,-0.0585) (-0.2020,-0.0586) (-0.1738,-0.1175) (-0.1738,-0.1175) (-0.1738,-0.1175) (-0.1617,-0.1338) (-0.1585,-0.1375) (-0.1585,-0.1375) (-0.1585,-0.1375) (-0.1024,-0.1829) (-0.1024,-0.1829) (-0.1024,-0.1829) (-0.0916,-0.1887) (-0.0832,-0.1924) (-0.0832,-0.1924) (-0.0832,-0.1924) (-0.0117,-0.2097) (-0.0116,-0.2097) (-0.0116,-0.2097) (-0.0048,-0.2099) (0.0109,-0.2096) (0.0109,-0.2096) (0.0110,-0.2096) (0.0789,-0.1952) (0.0789,-0.1952) (0.0789,-0.1952) (0.0822,-0.1938) (0.1044,-0.1827) (0.1044,-0.1827) (0.1044,-0.1827) (0.1525,-0.1458) (0.1527,-0.1456) (0.1527,-0.1456) (0.1538,-0.1443) (0.1760,-0.1158) (0.1760,-0.1158) (0.1760,-0.1158) (0.1983,-0.0721) (0.1986,-0.0714) (0.1986,-0.0714) (0.1988,-0.0707) (0.2093,-0.0233) (0.2093,-0.0233) (0.2093,-0.0233) (0.2102,0.0133) (0.2102,0.0151) (0.2102,0.0151) (0.2101,0.0153) (0.1967,0.0740) (0.1967,0.0740) (0.1967,0.0740) (0.1864,0.0971) (0.1847,0.1004) (0.1847,0.1004) (0.1846,0.1005) (0.1419,0.1546) (0.1419,0.1546) (0.1419,0.1546) (0.1297,0.1650) (0.1240,0.1693) (0.1240,0.1693) (0.1240,0.1693) (0.0586,0.2016) (0.0585,0.2016) (0.0585,0.2016) (0.0494,0.2041) (0.0371,0.2066) (0.0371,0.2066) (0.0370,0.2066) (-0.0344,0.2076) (-0.0344,0.2076) (-0.0344,0.2076) (-0.0395,0.2067) (-0.0594,0.2017) (-0.0594,0.2017) (-0.0594,0.2017) (-0.1183,0.1746) (-0.1185,0.1745) (-0.1185,0.1745) (-0.1205,0.1730) (-0.1443,0.1534) (-0.1443,0.1534) (-0.1443,0.1534) (-0.1792,0.1115) (-0.1796,0.1110) (-0.1796,0.1110) (-0.1796,0.1110) (-0.1801,0.1100) (-0.1983,0.0714) (-0.1983,0.0714) (-0.1983,0.0714) (-0.2087,0.0301) (-0.2090,0.0289) (-0.2090,0.0289) (-0.2090,0.0284) (-0.2088,-0.0263) (-0.2088,-0.0263) (-0.2088,-0.0263) (-0.2028,-0.0564) (-0.2021,-0.0591) (-0.2021,-0.0591) (-0.2020,-0.0593) (-0.1740,-0.1177) (-0.1740,-0.1177) (-0.1740,-0.1177) (-0.1618,-0.1341) (-0.1582,-0.1383) (-0.1582,-0.1383) (-0.1582,-0.1384) (-0.1029,-0.1831) (-0.1029,-0.1831) (-0.1029,-0.1831) (-0.0918,-0.1890) (-0.0826,-0.1931) (-0.0826,-0.1931) (-0.0826,-0.1931) (-0.0124,-0.2099) (-0.0124,-0.2099) (-0.0124,-0.2099) (-0.0052,-0.2102) (0.0116,-0.2099) (0.0116,-0.2099) (0.0116,-0.2099) (0.0781,-0.1958) (0.0782,-0.1957) (0.0782,-0.1957) (0.0817,-0.1942) (0.1049,-0.1826) (0.1049,-0.1826) (0.1049,-0.1826) (0.1521,-0.1464) (0.1524,-0.1461) (0.1524,-0.1461) (0.1524,-0.1461) (0.1536,-0.1448) (0.1763,-0.1155) (0.1763,-0.1155) (0.1763,-0.1155) (0.1983,-0.0725) (0.1986,-0.0716) (0.1986,-0.0716) (0.1986,-0.0716) (0.1988,-0.0708) (0.2094,-0.0230) (0.2094,-0.0230) (0.2094,-0.0230) (0.2104,0.0134) (0.2103,0.0154) (0.2103,0.0154) (0.2102,0.0157) (0.1968,0.0742) (0.1968,0.0743) (0.1968,0.0743) (0.1865,0.0974) (0.1846,0.1011) (0.1846,0.1011) (0.1845,0.1012) (0.1422,0.1547) (0.1422,0.1547) (0.1422,0.1547) (0.1298,0.1653) (0.1235,0.1701) (0.1235,0.1701) (0.1235,0.1701) (0.0592,0.2018) (0.0591,0.2018) (0.0591,0.2018) (0.0497,0.2043) (0.0364,0.2070) (0.0364,0.2070) (0.0364,0.2070) (-0.0336,0.2080) (-0.0337,0.2080) (-0.0337,0.2080) (-0.0390,0.2070) (-0.0600,0.2018) (-0.0600,0.2018) (-0.0600,0.2018) (-0.1177,0.1752) (-0.1179,0.1750) (-0.1179,0.1750) (-0.1179,0.1750) (-0.1202,0.1734) (-0.1447,0.1533) (-0.1447,0.1533) (-0.1447,0.1533) (-0.1790,0.1120) (-0.1794,0.1115) (-0.1794,0.1115) (-0.1794,0.1115) (-0.1801,0.1103) (-0.1985,0.0711) (-0.1985,0.0711) (-0.1985,0.0711) (-0.2088,0.0303) (-0.2091,0.0288) (-0.2091,0.0288) (-0.2091,0.0288) (-0.2091,0.0283) (-0.2089,-0.0266) (-0.2089,-0.0266) (-0.2089,-0.0266) (-0.2030,-0.0566) (-0.2021,-0.0596) (-0.2021,-0.0596) (-0.2021,-0.0599) (-0.1742,-0.1178) (-0.1742,-0.1179) (-0.1742,-0.1179) (-0.1619,-0.1344) (-0.1579,-0.1391) (-0.1579,-0.1391) (-0.1579,-0.1391) (-0.1034,-0.1831) (-0.1034,-0.1832) (-0.1034,-0.1832) (-0.0920,-0.1892) (-0.0821,-0.1936) (-0.0821,-0.1936) (-0.0820,-0.1937) (-0.0131,-0.2102) (-0.0131,-0.2102) (-0.0131,-0.2102) (-0.0056,-0.2105) (0.0122,-0.2101) (0.0122,-0.2101) (0.0122,-0.2101) (0.0774,-0.1962) (0.0776,-0.1962) (0.0776,-0.1962) (0.0776,-0.1962) (0.0813,-0.1946) (0.1054,-0.1825) (0.1054,-0.1825) (0.1054,-0.1825) (0.1517,-0.1469) (0.1521,-0.1466) (0.1521,-0.1466) (0.1521,-0.1466) (0.1534,-0.1451) (0.1766,-0.1153) (0.1766,-0.1153) (0.1766,-0.1153) (0.1982,-0.0728) (0.1987,-0.0718) (0.1987,-0.0718) (0.1987,-0.0718) (0.1989,-0.0709) (0.2096,-0.0227) (0.2096,-0.0227) (0.2096,-0.0227) (0.2105,0.0134) (0.2104,0.0157) (0.2104,0.0157) (0.2104,0.0157) (0.2104,0.0161) (0.1970,0.0744) (0.1970,0.0744) (0.1970,0.0744) (0.1866,0.0976) (0.1845,0.1018) (0.1845,0.1018) (0.1844,0.1019) (0.1426,0.1548) (0.1426,0.1548) (0.1426,0.1548) (0.1300,0.1656) (0.1231,0.1707) (0.1231,0.1707) (0.1230,0.1708) (0.0598,0.2019) (0.0597,0.2019) (0.0597,0.2019) (0.0500,0.2045) (0.0358,0.2074) (0.0358,0.2074) (0.0358,0.2074) (-0.0329,0.2083) (-0.0330,0.2083) (-0.0330,0.2083) (-0.0330,0.2083) (-0.0386,0.2073) (-0.0605,0.2018) (-0.0605,0.2018) (-0.0605,0.2018) (-0.1172,0.1757) (-0.1175,0.1755) (-0.1175,0.1755) (-0.1175,0.1755) (-0.1199,0.1738) (-0.1451,0.1531) (-0.1451,0.1531) (-0.1451,0.1531) (-0.1789,0.1125) (-0.1793,0.1118) (-0.1793,0.1118) (-0.1793,0.1118) (-0.1800,0.1106) (-0.1987,0.0708) (-0.1987,0.0708) (-0.1987,0.0708) (-0.2089,0.0305) (-0.2093,0.0288) (-0.2093,0.0288) (-0.2093,0.0288) (-0.2093,0.0281) (-0.2090,-0.0268) (-0.2090,-0.0268) (-0.2090,-0.0268) (-0.2031,-0.0567) (-0.2022,-0.0602) (-0.2022,-0.0602) (-0.2022,-0.0602) (-0.2021,-0.0604) (-0.1745,-0.1179) (-0.1744,-0.1179) (-0.1744,-0.1179) (-0.1620,-0.1346) (-0.1577,-0.1398) (-0.1577,-0.1398) (-0.1576,-0.1398) (-0.1039,-0.1832) (-0.1038,-0.1832) (-0.1038,-0.1832) (-0.0922,-0.1894) (-0.0815,-0.1942) (-0.0815,-0.1942) (-0.0815,-0.1942) (-0.0138,-0.2104) (-0.0137,-0.2104) (-0.0137,-0.2104) (-0.0137,-0.2104) (-0.0059,-0.2107) (0.0128,-0.2103) (0.0128,-0.2103) (0.0128,-0.2103) (0.0768,-0.1967) (0.0770,-0.1966) (0.0770,-0.1966) (0.0770,-0.1966) (0.0809,-0.1950) (0.1059,-0.1824) (0.1059,-0.1824) (0.1059,-0.1824) (0.1514,-0.1475) (0.1518,-0.1471) (0.1518,-0.1471) (0.1518,-0.1471) (0.1532,-0.1455) (0.1769,-0.1150) (0.1769,-0.1150) (0.1769,-0.1150) (0.1982,-0.0732) (0.1987,-0.0720) (0.1987,-0.0720) (0.1987,-0.0720) (0.1990,-0.0710) (0.2097,-0.0225) (0.2097,-0.0225) (0.2097,-0.0225) (0.2107,0.0133) (0.2106,0.0160) (0.2106,0.0160) (0.2106,0.0160) (0.2105,0.0165) (0.1971,0.0746) (0.1971,0.0746) (0.1971,0.0746) (0.1868,0.0978) (0.1844,0.1024) (0.1844,0.1024) (0.1844,0.1024) (0.1843,0.1026) (0.1429,0.1548) (0.1429,0.1548) (0.1429,0.1548) (0.1301,0.1657) (0.1227,0.1714) (0.1227,0.1714) (0.1226,0.1714) (0.0604,0.2020) (0.0603,0.2020) (0.0603,0.2020) (0.0503,0.2047) (0.0352,0.2078) (0.0352,0.2078) (0.0352,0.2078) (-0.0322,0.2086) (-0.0324,0.2086) (-0.0324,0.2086) (-0.0382,0.2075) (-0.0611,0.2019) (-0.0611,0.2019) (-0.0611,0.2019) (-0.1167,0.1762) (-0.1171,0.1760) (-0.1171,0.1760) (-0.1196,0.1741) (-0.1455,0.1529) (-0.1455,0.1529) (-0.1455,0.1528) (-0.1787,0.1129) (-0.1793,0.1121) (-0.1793,0.1121) (-0.1800,0.1108) (-0.1989,0.0705) (-0.1989,0.0705) (-0.1989,0.0705) (-0.2090,0.0307) (-0.2094,0.0287) (-0.2094,0.0287) (-0.2094,0.0280) (-0.2092,-0.0270) (-0.2092,-0.0270) (-0.2092,-0.0270) (-0.2033,-0.0567) (-0.2022,-0.0607) (-0.2022,-0.0607) (-0.2021,-0.0609) (-0.1747,-0.1180) (-0.1747,-0.1180) (-0.1747,-0.1180) (-0.1622,-0.1348) (-0.1574,-0.1404) (-0.1574,-0.1404) (-0.1573,-0.1405) (-0.1043,-0.1832) (-0.1043,-0.1833) (-0.1043,-0.1833) (-0.0924,-0.1895) (-0.0810,-0.1947) (-0.0809,-0.1947) (-0.0144,-0.2106) (-0.0143,-0.2106) (-0.0143,-0.2106) (-0.0063,-0.2109) (0.0134,-0.2105) (0.0134,-0.2105) (0.0762,-0.1971) (0.0765,-0.1970) (0.0765,-0.1970) (0.0806,-0.1952) (0.1064,-0.1823) (0.1064,-0.1823) (0.1510,-0.1480) (0.1516,-0.1475) (0.1516,-0.1475) (0.1531,-0.1457) (0.1772,-0.1148) (0.1772,-0.1148) (0.1982,-0.0735) (0.1988,-0.0721) (0.1988,-0.0721) (0.1991,-0.0710) (0.2099,-0.0222) (0.2099,-0.0222) (0.2108,0.0132) (0.2107,0.0163) (0.2107,0.0163) (0.2106,0.0168) (0.1973,0.0746) (0.1973,0.0746) (0.1870,0.0979) (0.1843,0.1030) (0.1843,0.1030) (0.1841,0.1032) (0.1433,0.1548) (0.1433,0.1548) (0.1303,0.1658) (0.1223,0.1719) (0.1223,0.1719) (0.1222,0.1719) (0.0609,0.2020) (0.0609,0.2020) (0.0507,0.2048) (0.0347,0.2081) (0.0347,0.2081) (-0.0315,0.2089) (-0.0315,0.2089) (-0.0376,0.2077)};
\addplot[black,dashed,domain=0:360,samples=181] ({0.2158*cos(x)},{0.2158*sin(x)});
\addplot[only marks,mark=square*,mark size=1.6pt,blue!70!black] coordinates {(-0.0257,-0.0000) (0.0060,-0.0000) (-0.0014,0.0000) (-0.0001,0.0000) (0.0003,0.0000)};
\node[red!80!black,font=\scriptsize,anchor=north] at (axis cs:0.0789,-0.1952) {$T=50$};
\addplot[only marks,mark=*,mark size=1.5pt,red!80!black] coordinates {(0.1075,-0.1218)};
\node[blue!70!black,font=\scriptsize,anchor=south,align=center] at (axis cs:-0.01,0.025) {trace ordering};
\end{axis}
\end{tikzpicture}
\caption{(a) The values $\frac12\tr\rho([C_\gamma])$ over the $2{,}203{,}681$ modular knots with $\tr\gamma\leq 6000$, in $50$ bins,
for the homomorphism $\rho$ of \Cref{rem.numerics.su2} with $\tr\rho(xy)=\frac12$,
together with the density $\frac{2}{\pi}\sqrt{1-x^{2}}$ of \Cref{cor.satotate}.
(b) The averages of $\rho_\theta([C_\gamma])$ of \Cref{eg.circle.failure} in the complex plane, with $\theta=(1+\varepsilon)/2$ and $\varepsilon=0.05\sqrt2$.
The red path is the average over the classes with $\ell_h(\gamma)\leq T$ as $T$ increases from $10$ to $59$,
and the dashed circle has radius $\delta/\sqrt{\delta^{2}+\pi^{2}}=0.2158$.
The red path starts at the dot, where $T=10$.
The blue squares are the averages over the classes with $\tr\gamma\leq N$ for $N=60$, $200$, $600$, $2000$, and $6000$; the last three overlap.}
\label{fig.numerics}
\end{figure}

\begin{rem}[Numerical illustration] \label{rem.numerics.circle}
For $\varepsilon=0.05\sqrt2$ the entropy is $\delta\approx 0.6944$ and the limit is $0.2158$;
the moduli computed by counting the primitive classes with M\"obius inversion are $0.2056$, $0.2090$, and $0.2106$ for $T=30$, $40$,
and $50$, over $5.9\times 10^{7}$, $4.5\times 10^{10}$, and $3.7\times 10^{13}$ classes,
and they increase towards the limit with a correction of order $1/T$,
in accordance with the next term of the asymptotic expansion of the sum.
In contrast, the logarithmic density of \Cref{thmA.compact} \textup{(4)} converges more slowly:
with the modular knots indexed from $i=1$, the sums $\sum_{i\leq N}\rho_\theta([K_i])/i$ at the same three thresholds have moduli $0.415$, $0.409$,
and $0.404$, in accordance with the bound $O(1)$ in the proof of \Cref{thmA.general} \textup{(2)},
so the averages $\frac{1}{\log N}\sum_{i\leq N}\rho_\theta([K_i])/i$ decrease like $1/\log N$.
The constant depends on the indexing: when $K_{2,3}$ occupies the first position,
the moduli are $0.183$, $0.186$, and $0.187$.
For a character outside \Cref{prop.exceptional.explicit}, for instance $\theta=\sqrt3-1$,
the natural averages of $\rho_\theta$ over $\ell_h\leq T$ are $0.046$, $0.044$, and $0.042$, decreasing to $0$.
For the trace ordering, on the other hand,
the average of $\rho_\theta([C_\gamma])$ with $\theta=(1+\varepsilon)/2$ over the classes with $\tr\gamma\leq N$ is $0.026$, $0.006$,
and $0.0014$ for $N=60$, $200$, and $600$, in accordance with \Cref{thmA.compact} \textup{(1)}.
\Cref{fig.numerics} (b) shows the two kinds of averages.
\end{rem}

\begin{rem} \label{rem.mcmullen.circle}
McMullen \cite[Section 6]{McMullen2013CM} observed that the law for $G=S^{1}$ can fail,
the example being a character that is a function of the period of a suspension flow.
\Cref{eg.circle.failure} shows that the failure occurs for the modular knots and a generic ordering, with a character of infinite image.
\Cref{thmA.compact} \textup{(1)} shows that it does not occur for the trace ordering:
the cusp, which makes the hyperbolic length of $R^{q}L$ logarithmic in $q$,
prevents any character from being a function of the hyperbolic length.
\end{rem}

\begin{eg} \label{eg.circle.su2}
Let $G=S^{1}\times \mathrm{SU}(2)$ and $\rho=(\rho_\theta,\rho')\colon \pi_1(S^{3}-K_{2,3})\to G$ with $\theta$ as in \Cref{eg.circle.failure} and $\rho'$ with dense image in $\mathrm{SU}(2)$;
the image of $\rho$ is dense in $G$, because the closure of the image projects onto both factors and $\mathrm{SU}(2)$ has no nontrivial abelian quotient.
The special characters are the $\chi^{k}$ with $\chi$ the projection to $S^{1}$,
so $G_S=\mathrm{SU}(2)$, $A=S^{1}$, and $\eta(t)=e^{i\pi t}$.
\Cref{prop.exceptional.fibers} says that, for the ordering by $\ell_h=\ell_{\rm c}+\varepsilon\Psi$,
the pairs $\bigl(e^{2\pi i\theta\Psi(\gamma)},\tfrac12\tr\rho'([C_\gamma])\bigr)$ are distributed as follows:
the first coordinate equals $e^{i\pi\ell_h(\gamma)}$ exactly,
while the second is equidistributed with respect to $\frac{2}{\pi}\sqrt{1-x^{2}}\,dx$,
jointly with any continuous function of the first coordinate.
The Sato--Tate law for the $\mathrm{SU}(2)$ component survives the failure of the law for the $S^{1}$ component.
\end{eg}

\begin{rem}[Numerical illustration] \label{rem.numerics.su2}
Take $\rho'$ with $\rho'(x)=\mathrm{diag}(i,-i)$ and with $\rho'(y)$ such that $\rho'(y)^{3}=-I$, $\tr\rho'(y)=1$, and $\tr\rho'(xy)=c$,
and compute $\rho'([C_\gamma])$ by substituting $\rho'(xy^{-1})$ for $R$ and $\rho'(xy^{-2})$ for $L$ in the word of $\gamma$, as in \Cref{eg.su2.brieskorn}.
We take $c=2\cos\frac{2\pi}{7}$, the value for the homomorphism $\rho_1$ of \Cref{eg.su2.brieskorn},
and $c=\frac12$; in both cases $\rho'$ has dense image, by \Cref{sec.su2}.
For each pair $(n_L,n_R)$, the sum of $\frac12\tr\rho'([C_\gamma])$ over the classes with $n_L(\gamma)=n_L$
and $n_R(\gamma)=n_R$ is computed exactly from the sums of the traces of the words,
by M\"obius inversion over the non-primitive words,
where $\tr(g^{k})=\tr\operatorname{Sym}^{k}g -\tr\operatorname{Sym}^{k-2}g$ for $g\in \mathrm{SU}(2)$ is used for the powers.
With $\varepsilon=0.05\sqrt2$ and $T=26$, $34$, $42$, there are $4.3\times10^{6}$, $8.3\times10^{8}$,
and $1.7\times10^{11}$ classes with $\ell_h(\gamma)\leq T$,
and the averages of $e^{i\pi\ell_h(\gamma)}$ have moduli $0.203$, $0.207$, $0.209$.
The averages of the products $e^{i\pi\ell_h(\gamma)}\,\frac12\tr\rho'([C_\gamma])$ have moduli $0.0075$, $0.0020$, $0.00055$ for $c=2\cos\frac{2\pi}{7}$ and $0.0033$, $0.0011$, $0.00046$ for $c=\frac12$,
decreasing to $0$ as \Cref{prop.exceptional.fibers} predicts,
whereas a correlation between the two coordinates would keep them of the size of the first averages.
\end{rem}

\subsection{The combinatorial ordering} \label{subsec.compact.comb}

\begin{proof}[Proof of \Cref{thmA.compact} \textup{(3)}]
The proof of \PI{thm.combinatorial} applies to each irreducible unitary representation $\varrho$ of $G$ separately, as follows.
Let $\mca{G}'$ be the refined graph, and let $L_\varrho$ denote the adjacency matrix twisted by $\varrho\circ\tilde\rho$, as in \PI{lem.twisted}.
That lemma assumes that the holonomies of the directed loops at a vertex generate $G$;
here they generate a dense subgroup, by \PI{lem.dense} \textup{(1)},
and that suffices for its proof, in which the generation is used only to see that a line fixed by all holonomies is fixed by $G$, because $\varrho$ is continuous.
Hence the spectral radius of $L_\varrho$ is at most $2$,
with equality if and only if $\varrho$ is a character $\chi$ and there is $\lambda\in S^{1}$ with $\chi([\tau])=\lambda^{n}$ for every periodic orbit $\tau$ of period $n$.
Evaluating along the closed paths of the classes $R^{a}L$ and $RL^{b}$ as in Step 3 of the proof of \PI{thm.combinatorial} yields $\lambda=\pm 1$,
and then $\chi(\rho([C_\gamma]))=\lambda^{\ell_{\rm c}(\gamma)} =\lambda^{\operatorname{lk}(C_\gamma,K_{2,3})}$ for all but finitely many $\gamma$, because $\ell_{\rm c}\equiv \Psi\pmod 2$.

If $\lambda=1$, then $\chi\circ\rho$ is trivial on all but finitely many modular knots,
hence $\chi=1$ by \Cref{lem.generate.linking} and density.
If $\lambda=-1$, then $\chi^{2}\circ\rho$ is trivial on all but finitely many modular knots,
hence on $\pi_1(S^{3}-L_n)$ by \Cref{lem.generate.linking},
so that $\chi^{2}=1$ on the dense image of $\rho$,
hence on $G$, and $\chi$ is a character as in the statement.
Conversely, such a $\chi$ satisfies $\chi(\rho([C_\gamma]))=(-1)^{n}$ for $\ell_{\rm c}(\gamma)=n$,
so that its average over $\ell_{\rm c}\leq N$ oscillates between $\frac13$
and $-\frac13$ by the proof of \PI{thm.combinatorial}, and the law fails for $F=\chi$.

Let $\varrho$ be an irreducible unitary representation for which the spectral radius of $L_\varrho$ is smaller than $2$, so that the trace of $L_\varrho^{n}$ is $O(r_0^{n})$ with $r_0<2$.
By \PIeq{eq.trace}, that trace is the sum of $\chi_\varrho([x])$ over the periodic points $x$ with $\sigma^{n}(x)=x$,
and the periodic points whose least period is smaller than $n$ contribute $O(2^{n/2})$ to it, as in Step 2 of that proof.
Hence the sum of $\chi_\varrho(\rho([C_\gamma]))$ over the modular knots with $\ell_{\rm c}(\gamma)=n$ is $O(r_0^{n}/n)+O(2^{n/2})$, which is $o(2^{n}/n)$.
If no such $\chi$ exists, then this applies to every nontrivial $\varrho$,
and the argument of Step 4 of that proof yields the equality for $F$ the character of $\varrho$;
the equality for every continuous class function follows by the Peter--Weyl theorem.
\end{proof}

For $G=\mathrm{SU}(2)$ there is no nontrivial character,
so the modular knots ordered by the combinatorial length,
with the missing trefoil inserted or not, satisfy the conclusion of \Cref{thmA.compact} \textup{(3)} for every $\rho$ with dense image.

\begin{rem} \label{rem.exceptional.comb}
The combinatorial ordering has the same structure.
A character $\chi\colon G\to\{\pm1\}$ with $\chi\circ\rho=(-1)^{\operatorname{lk}(\cdot,K_{2,3})}$ is unique if it exists,
because the quotient of two such characters is trivial on all but finitely many modular knots,
hence trivial by \Cref{lem.generate.linking} and density.
The proof of \Cref{thmA.compact} \textup{(3)} shows,
for every irreducible $\varrho$ other than $1$ and that $\chi$,
that the sum of $\chi_\varrho$ over the modular knots with $\ell_{\rm c}=n$ is $o(2^{n}/n)$, whether $\chi$ exists or not.
Hence, if $\chi$ exists, the Frobenius classes of the modular knots with $\ell_{\rm c}=n$ are equidistributed in the fiber of $G\to\{\pm1\}$ over $(-1)^{n}$,
in the sense of \Cref{prop.exceptional.fibers} with the level $n$ in place of the threshold $T$.
\end{rem}

\subsection{Logarithmic density for the three orderings} \label{subsec.compact.log.modular}

\begin{proof}[Proof of \Cref{thmA.compact} \textup{(4)}]
For the trace ordering, \eqref{eq.law.compact} holds by \textup{(1)},
and \Cref{lem.log.partial} implies \eqref{eq.law.log}.
For a generic ordering, this is \Cref{thmA.general} \textup{(2)}.
For the combinatorial ordering, let $\varrho$ be a nontrivial irreducible unitary representation.
If $\varrho$ is not a character with $\chi\circ\rho=(-1)^{\operatorname{lk}(\cdot,K_{2,3})}$,
then the average of $\chi_\varrho$ over the level $\ell_{\rm c}=m$ tends to $0$,
by \Cref{rem.exceptional.comb}; otherwise that average is $(-1)^{m}+o(1)$, whose Ces\`aro means tend to $0$.
The assertion follows by the Peter--Weyl theorem.
\end{proof}

The combinatorial ordering is stated in the form of level averages rather than in the form \eqref{eq.law.log} because the knots of a given combinatorial length are about half of all the knots up to that length,
so that a sequence form depends on the order chosen inside a level,
as noted in \PI{rem.combinatorial.ties}; the level averages do not.
Since $P_m\sim 2^{m}/m$, the number $m$ is $\log_2 N$ up to an error $O(\log\log N)$,
and the average over $m\leq m_0$ is an average of logarithmic type.

\section{Representations into $\mathrm{SU}(2)$} \label{sec.su2}

We first describe a coordinate on the space of the representations of \Cref{cor.satotate},
and we then locate on it two families that arise in other contexts.
The conjugacy classes of the homomorphisms $\rho:\pi_1(S^{3}-K_{2,3})\to \mathrm{SU}(2)$ with nonabelian image form a single open arc,
by Klassen \cite[Theorem 1]{Klassen1991TransAMS}.
Write $\pi_1(S^{3}-K_{2,3})=\langle x,y\mid x^{2}=y^{3}\rangle$,
let $\rho$ be such a homomorphism, put $(A,B)=(\rho(x),\rho(y))$,
and let $\ol\rho$ be the composition of $\rho$ with the projection $\mathrm{SU}(2)\to \mathrm{SO}(3)$.
Then $A^{2}=B^{3}=-I$, so that $\tr A=0$ and $\tr B=1$, and the class of $\rho$ is determined by $\tr(AB)$.
That value is a coordinate on the arc, with values in $(-\sqrt3,\sqrt3)$.
The abelianization of $\pi_1(S^{3}-K_{2,3})$ is infinite cyclic,
so there is a unique nontrivial homomorphism $\epsilon\colon \pi_1(S^{3}-K_{2,3})\to\{\pm1\}$, the sign character.
It sends $x$ to $-1$ and $y$ to $1$, so that $\tr(\epsilon\rho)(xy)=-\tr\rho(xy)$,
and $\epsilon([C_\gamma])=(-1)^{\Psi(\gamma)}$ by \PI{thm.ghys}.

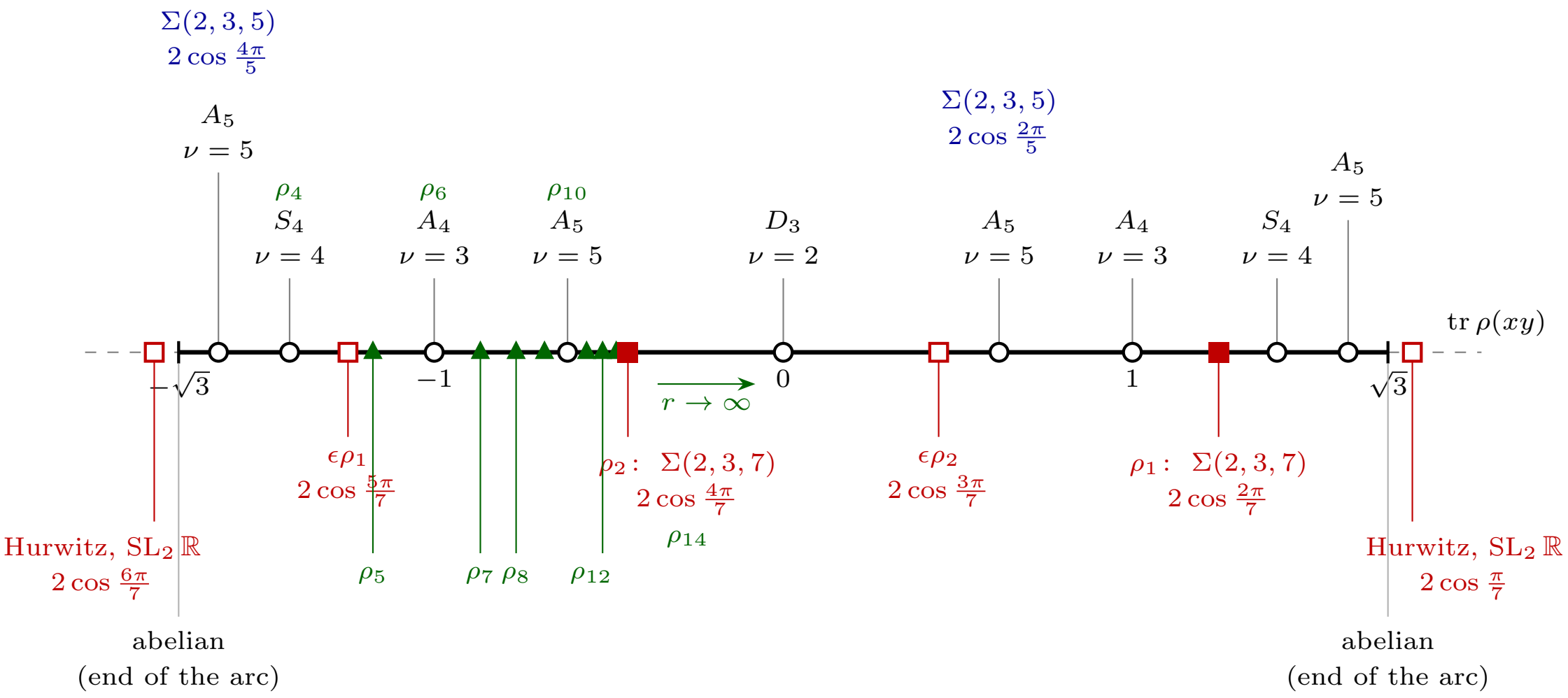
\begin{figure}[t]
\centering
\begin{tikzpicture}[x=3.3cm,y=1cm,>=Stealth,font=\fontsize{7.5}{9.5}\selectfont]
  \def\sq{1.7320508}
  \draw[dashed,gray] (-2.0,0) -- (-\sq,0);  \draw[dashed,gray] (\sq,0) -- (2.0,0);
  \draw[very thick] (-\sq,0) -- (\sq,0);
  \node[above right] at (1.86,0.02) {$\tr\rho(xy)$};
  \foreach \t/\l in {-1/{-1},0/{0},1/{1}} { \draw (\t,2pt) -- (\t,-2pt); \node[below=1pt] at (\t,0) {$\l$}; }
  \foreach \t/\l in {-\sq/{-\sqrt3},\sq/{\sqrt3}} {
    \draw[thick] (\t,3pt) -- (\t,-3pt); \node[below=1pt] at (\t,0) {$\l$};
    \draw[gray!60] (\t,-3pt) -- (\t,-2.5); \node[below,align=center] at (\t,-2.5) {abelian\\(end of the arc)}; }
  \foreach \t/\g/\n/\h in {0/{D_3}/2/0.7, -1/{A_4}/3/0.7, 1/{A_4}/3/0.7, -1.4142136/{S_4}/4/0.7, 1.4142136/{S_4}/4/0.7,
                        -0.6180340/{A_5}/5/0.7, 0.6180340/{A_5}/5/0.7, -1.6180340/{A_5}/5/1.7, 1.6180340/{A_5}/5/1.25} {
    \draw[gray] (\t,0) -- (\t,\h);
    \draw[thick,fill=white] (\t,0) circle[radius=2.4pt];
    \node[above,align=center] at (\t,\h) {$\g$\\$\nu=\n$}; }
  \node[blue!60!black,above,align=center] at (0.6180340,1.75) {$\Sigma(2,3,5)$\\$2\cos\frac{2\pi}{5}$};
  \node[blue!60!black,above,align=center] at (-1.6180340,2.5) {$\Sigma(2,3,5)$\\$2\cos\frac{4\pi}{5}$};
  \node[green!40!black,above] at (-1.4142136,1.3) {$\rho_4$};
  \node[green!40!black,above] at (-1,1.3) {$\rho_6$};
  \node[green!40!black,above] at (-0.6180340,1.3) {$\rho_{10}$};
  \foreach \r/\t in {5/-1.1755705, 7/-0.8677675, 8/-0.7653669, 9/-0.6840403, 11/-0.5634651, 12/-0.5176381, 13/-0.4786313} {
    \fill[green!40!black] (\t,0) ++(0,3.2pt) -- ++(-2.6pt,-5pt) -- ++(5.2pt,0) -- cycle; }
  \foreach \t/\l/\nm/\dx in {1.2469796/{2\cos\frac{2\pi}{7}}/{\rho_1\colon\ \Sigma(2,3,7)}/0pt, -0.4450419/{2\cos\frac{4\pi}{7}}/{\rho_2\colon\ \Sigma(2,3,7)}/16pt} {
    \draw[red!75!black] (\t,0) -- (\t,-0.8);
    \fill[red!75!black] (\t,0) ++(-2.8pt,-2.8pt) rectangle ++(5.6pt,5.6pt);
    \node[red!75!black,below,align=center,xshift=\dx] at (\t,-0.8) {$\nm$\\$\l$}; }
  \node[green!40!black,below,xshift=16pt] at (-0.4450419,-1.55) {$\rho_{14}$};
  \foreach \t/\l/\nm in {-1.2469796/{2\cos\frac{5\pi}{7}}/{\epsilon\rho_1}, 0.4450419/{2\cos\frac{3\pi}{7}}/{\epsilon\rho_2}} {
    \draw[red!75!black] (\t,0) -- (\t,-0.8);
    \draw[red!75!black,thick,fill=white] (\t,0) ++(-2.4pt,-2.4pt) rectangle ++(4.8pt,4.8pt);
    \node[red!75!black,below,align=center] at (\t,-0.8) {$\nm$\\$\l$}; }
  \foreach \t/\l/\dx in {-1.8019377/{2\cos\frac{6\pi}{7}}/-14pt,1.8019377/{2\cos\frac{\pi}{7}}/14pt} {
    \draw[red!75!black] (\t,0) -- (\t,-1.6);
    \draw[red!75!black,thick,fill=white] (\t,0) ++(-2.4pt,-2.4pt) rectangle ++(4.8pt,4.8pt);
    \node[red!75!black,below,align=center,xshift=\dx] at (\t,-1.6) {Hurwitz, $\SL_2\R$\\$\l$}; }
  \foreach \r/\t in {5/-1.1755705, 7/-0.8677675, 8/-0.7653669} {
    \draw[green!40!black] (\t,0) -- (\t,-1.9); \node[green!40!black,below] at (\t,-1.9) {$\rho_{\r}$}; }
  \draw[green!40!black] (-0.5176381,0) -- (-0.5176381,-1.9); \node[green!40!black,below] at (-0.55,-1.9) {$\rho_{12}$};
  \draw[->,green!40!black] (-0.36,-0.3) -- (-0.08,-0.3) node[pos=0.5,below] {$r\to\infty$};
\end{tikzpicture}
\caption{The arc of the classes of homomorphisms $\rho\colon\pi_1(S^{3}-K_{2,3})\to\mathrm{SU}(2)$ with nonabelian image, in the coordinate $\tr\rho(xy)\in(-\sqrt3,\sqrt3)$.
The image of $\rho$ is dense except at the nine points marked by circles, where the image of $\ol\rho$ is the dihedral group of order six or a polyhedral group;
two of them are the classes through $\pi_1(\Sigma(2,3,5))$.
The filled squares are the two classes $\rho_1$ and $\rho_2$ through $\pi_1(\Sigma(2,3,7))$ of \Cref{eg.su2.brieskorn},
the hollow squares on the arc are their twists by the sign character $\epsilon$,
and the hollow squares outside the arc are the third conjugates, which are Fuchsian.
The triangles are the Jones representations $\rho_r$ of \Cref{eg.su2.jones}, which accumulate at $0$;
$\rho_4$, $\rho_6$, and $\rho_{10}$ are three of the nine exceptional points, and $\rho_{14}=\rho_2$.}
\label{fig.arc}
\end{figure}

Let $\nu\in \Z_{\geq 2}\cup\{\infty\}$ be the order of $\ol\rho(xy)$.
The image of $\ol\rho$ is generated by a rotation of order two and a rotation of order three about a distinct axis,
so it is a quotient of the triangle group $\langle a,b\mid a^{2}=b^{3}=(ab)^{\nu}=1\rangle$,
and $|\tr(AB)|=2\cos(k\pi/\nu)$ for some $k\in \Z_{\geq 1}$ with $\gcd(k,\nu)=1$.
By the classification of the closed subgroups of $\mathrm{SO}(3)$,
the subgroup generated by two such rotations is dense in $\mathrm{SO}(3)$ unless it is the dihedral group of order six or the rotation group of the regular tetrahedron, octahedron,
or icosahedron, and in those four groups the product of the two rotations has order two, three, four, and five respectively.
Since $\mathrm{SU}(2)$ is connected, it has no proper closed subgroup of finite index,
so the image of $\rho$ is dense in $\mathrm{SU}(2)$ if and only if the image of $\ol\rho$ is dense in $\mathrm{SO}(3)$.
The image of $\rho$ is therefore dense except at the nine points of the arc with
\[\tr(AB)\in\Bigl\{0,\ \pm1,\ \pm\sqrt2,\ \pm\tfrac{\sqrt5-1}{2},\
\pm\tfrac{\sqrt5+1}{2}\Bigr\}.\]
In particular, for $\nu\in \Z_{\geq 7}$ and $k\in \Z_{\geq 1}$ with $\gcd(k,\nu)=1$ and $|2\cos(k\pi/\nu)|<\sqrt3$,
the representations in the class with $\tr(AB)=2\cos(k\pi/\nu)$ have dense image.
\Cref{fig.arc} shows the arc, the nine exceptional points, and the two families described below.

The first family comes from the Brieskorn spheres.
A meridian of $K_{2,3}$ is $\mu=xy^{-1}$ and a longitude is $\lambda=x^{2}\mu^{-6}$,
as for the $(3,4)$-torus knot in the proof of Klassen \cite[Proposition 2]{Klassen1991TransAMS}.
Let $n\in \Z$, and let $\Sigma$ be the homology sphere obtained by the Dehn surgery along $\mu\lambda^{n}$.
Since $\rho(x^{2})=-I$, the representation $\rho$ factors through $\pi_1(\Sigma)$ if and only if $\rho(\mu)^{1-6n}=(-1)^{n}I$,
and $\tr\rho(\mu)=-\tr(AB)$ by the identity $\tr(AB)+\tr(AB^{-1})=\tr A\,\tr B$.
For $n=1$ the sphere $\Sigma$ is the Poincar\'e sphere $\Sigma(2,3,5)$,
and the two classes that satisfy the condition are those with $\tr(AB)=2\cos(2\pi/5)$ and $\tr(AB)=2\cos(4\pi/5)$, the two roots of $X^{2}+X-1$;
they are two of the nine exceptional points above,
at which the image of $\rho$ is the binary icosahedral group of order $120$.
The triangle group with $\nu=5$ is spherical,
and accordingly the field $\Q(\cos(2\pi/5))$ is quadratic and neither of the two roots lies outside $[-\sqrt3,\sqrt3]$.

\begin{eg} \label{eg.su2.brieskorn}
Take $n=-1$, so that $\Sigma=\Sigma(2,3,7)$.
There are four classes with $\nu=7$, given by $k\in\{2,3,4,5\}$,
and exactly two of them satisfy the condition above, namely those with
\[\tr(AB)=2\cos(2\pi/7)\qquad\text{and}\qquad \tr(AB)=2\cos(4\pi/7).\]
Their number agrees with $-2$ times the Casson invariant of $\Sigma(2,3,7)$, which is $-1$,
as it must by Fintushel and Stern \cite[Theorem 2.10]{FintushelStern1990ProcLMS};
they parametrize the conjugacy classes of the irreducible representations of the fundamental group of a Brieskorn sphere into $\mathrm{SU}(2)$ by triples of integers in \cite[Proposition 2.8]{FintushelStern1990ProcLMS},
and they count them in \cite[Corollary 2.9]{FintushelStern1990ProcLMS}.

Put $c=2\cos(2\pi/7)$, and let $\rho_1$ and $\rho_2$ lie in the two classes above, in that order.
The minimal polynomial of $c$ is $X^{3}+X^{2}-2X-1$,
whose three roots are $2\cos(2\pi/7)$, $2\cos(4\pi/7)$, and $2\cos(6\pi/7)$.
The two classes are given by the first two roots.
The third root has absolute value greater than $\sqrt3$,
and the corresponding pair $(A,B)$ lies in $\mathrm{SL}_2\R$ and defines,
up to the sign character, the Fuchsian representation onto the Hurwitz triangle group.

Every trace $\tr\rho_1([C_\gamma])$ is a $\Z$-polynomial in $\tr A$, $\tr B$, and $c$, hence an element of $\Z[c]$,
and $\tr\rho_2([C_\gamma])$ is its conjugate under $c\mapsto 2\cos(4\pi/7)$.
The class $[C_\gamma]$ is the unique lift of $\gamma$ whose image in $H_1(S^{3}-K_{2,3})\cong \Z$ is $\Psi(\gamma)$,
by \PI{thm.ghys}, so that $R$ lifts to $\mu$ and $L$ lifts to $x y^{-2}$.
Substituting $-AB^{2}$ for $R$ and $-AB$ for $L$ in the word of $\gamma$, we obtain the following values.
\begin{center}
\begin{tabular}{l|l|l|l|l|l}
$\gamma$ & $\tr\gamma$ & $\tr\rho_1([C_\gamma])$
 & at $c$ & at $2\cos\frac{4\pi}{7}$ & at $2\cos\frac{6\pi}{7}$ \\[2pt]
 \hline 
$LR$ & $3$ & $c^{2}-1$ & $0.5550$ & $-0.8019$ & $2.2470$\\
$L^{3}R$ & $5$ & $-c$ & $-1.2470$ & $0.4450$ & $1.8019$\\
$L^{2}R^{2}$ & $6$ & $1-c$ & $-0.2470$ & $1.4450$ & $2.8019$\\
$L^{4}R$ & $6$ & $1$ & $1.0000$ & $1.0000$ & $1.0000$\\
$L^{2}RLR$ & $10$ & $c^{2}$ & $1.5550$ & $0.1981$ & $3.2470$\\
$L^{2}RLR^{2}$ & $15$ & $c^{2}-c-1$ & $-0.6920$ & $-0.3569$ & $4.0489$
\end{tabular}
\end{center}
The first two of the three conjugates lie in $[-2,2]$,
and \Cref{cor.satotate} applies to $\rho_1$ and to $\rho_2$, so that each of the two columns equidistributes.
The third conjugate is the trace of the image of $\gamma$ in the Hurwitz group.
\end{eg}

The second family comes from the Jones representations.
Let $r\in \Z_{\geq 3}$, put $\zeta=e^{2\pi i/r}$,
and let $\rho_r$ be the two-dimensional sector of the Jones representation of $B_3=\pi_1(S^{3}-K_{2,3})$ at $\zeta$,
that is, the sector of the Burau hook, normalized so as to take values in $\mathrm{SU}(2)$.
Then $\tr\rho_r(xy)=-2\sin(\pi/r)$, so that the order $\nu$ above equals $2r$ if $r$ is odd,
equals $r/2$ if $r\equiv 2\pmod 4$, and equals $r$ if $r\equiv 0\pmod 4$.
The image of $\rho_3$ is abelian, and $\nu\leq 5$ exactly for $r\in\{4,6,10\}$.
The image of $\rho_r$ is therefore dense for every other value of $r$,
and $\{3,4,6,10\}$ is the list of exceptions in the density theorem of Freedman, Larsen,
and Wang \cite[Theorem 4.1]{FreedmanLarsenWang2002CMP}.

\begin{eg} \label{eg.su2.jones}
Let $r\in\Z_{\geq 3}$, let $\zeta=e^{2\pi i/r}$, and let $k\in \Z_{\geq 1}$ satisfy $\gcd(k,r)=1$.
The Burau matrices at $\zeta^{k}$ satisfy $A^{2}=B^{3}=-I$ and $\tr(AB)=-2\sin(k\pi/r)$,
and an irreducible representation is determined up to conjugation by its character,
so those matrices are conjugate into $\mathrm{SU}(2)$ if and only if $|2\sin(k\pi/r)|<\sqrt3$.
For $r=14$ the values $k=1$ and $k=3$ yield $\tr(AB)=-2\cos(3\pi/7)$ and $\tr(AB)=-2\cos(2\pi/7)$,
which are the two classes of \Cref{eg.su2.brieskorn} up to the sign character $\epsilon$, that is, the two classes of $\ol\rho$.
The value $k=5$ yields $\tr(AB)=-2\cos(\pi/7)$, which is the Fuchsian representation.
The three primitive fourteenth roots of unity, taken up to complex conjugation,
therefore correspond to the three real places of $\Q(\cos(\pi/7))$.
\end{eg}

Freedman, Larsen, and Wang also determine the limiting distribution of the evaluations of the Jones polynomial as the length of the braid grows \cite[Section 5]{FreedmanLarsenWang2002CMP}.
\Cref{cor.satotate} concerns instead the modular knots ordered by $|\tr\gamma|$.

\begin{proof}[Proof of \Cref{thmA.joint}]
Put $G=\mathrm{SU}(2)^{j}\times \Z/m\Z$, whose identity component $\mathrm{SU}(2)^{j}$ is semisimple,
and let $\rho$ be the product of the $\rho_i$ with the reduction modulo $m$ of the abelianization.
That reduction is surjective, because $\pi_1(S^{3}-K_{2,3})^{\rm ab}\cong \Z$,
and the last coordinate of $\rho([C_\gamma])$ is $\Psi(\gamma)\bmod m$, by \PI{thm.ghys}.
Let $H$ be the closure of the image of $\rho$, so that $H$ surjects onto every factor of $G$.

The closed normal subgroups of $\mathrm{SU}(2)$ are the trivial group, $\{\pm I\}$,
and $\mathrm{SU}(2)$, and every automorphism of $\mathrm{SU}(2)$ and of $\mathrm{SO}(3)$ is inner.
By Goursat's lemma, the closure of the image of $(\rho_i,\rho_{i'})$ in $\mathrm{SU}(2)^{2}$ is therefore either the whole group,
or the graph of an inner automorphism of $\mathrm{SU}(2)$,
or the preimage of the graph of an inner automorphism of $\mathrm{SO}(3)$.
In the last two cases $\rho_{i'}$ is conjugate to $\rho_i$ or to $\epsilon\rho_i$,
so that $|\tr\rho_{i'}(xy)|=|\tr\rho_i(xy)|$, which the assumption excludes.
We show by induction on $j$ that the closure $H_j$ of the image of $(\rho_1,\dots,\rho_j)$ is $\mathrm{SU}(2)^{j}$.
Suppose that $H_{j-1}=\mathrm{SU}(2)^{j-1}$.
The group $H_j$ projects onto $\mathrm{SU}(2)^{j-1}$ and onto the last factor,
so by Goursat's lemma there are closed normal subgroups $N_1\subset\mathrm{SU}(2)^{j-1}$ and $N_2\subset\mathrm{SU}(2)$ and an isomorphism $\mathrm{SU}(2)^{j-1}/N_1\cong\mathrm{SU}(2)/N_2$ whose graph is the image of $H_j$.
If $N_2\neq\mathrm{SU}(2)$, then the quotient is of dimension three, so that $N_1$ contains the product of all the factors of $\mathrm{SU}(2)^{j-1}$ except one, say the $i$-th one;
then the closure of the image of $(\rho_i,\rho_j)$ is not $\mathrm{SU}(2)^{2}$, which contradicts the case of two factors.
Hence $N_2=\mathrm{SU}(2)$, and $H_j=\mathrm{SU}(2)^{j}$.

Let $N=\{h\in \mathrm{SU}(2)^{j}\mid (h,0)\in H\}$, a closed normal subgroup of $\mathrm{SU}(2)^{j}$.
The quotient $\mathrm{SU}(2)^{j}/N$ embeds into $\Z/m\Z$ and is therefore finite,
and $\mathrm{SU}(2)^{j}$ is connected, so $N=\mathrm{SU}(2)^{j}$ and $H=G$.
\Cref{thmA.compact} \textup{(1)} applies to the trace ordering.
For a generic ordering by $\ell_h$, every continuous character of $G$ has finite image,
because the identity component of $G$ is semisimple,
so that \eqref{eq.exceptional} cannot hold with $a\neq 0$,
as noted at the beginning of \Cref{sec.criterion}, nor with $a=0$ for a nontrivial character,
by \Cref{thmA.general} \textup{(1)}; hence \Cref{thmA.compact} \textup{(2)} applies.
The push-forward of the Haar measure of $G$ is the stated product measure, by the Weyl integration formula.
\end{proof}

\begin{rem} \label{rem.joint.assumption}
The assumption of \Cref{thmA.joint} on the values $|\tr\rho_i(xy)|$ cannot be dropped.
If $|\tr\rho_{i'}(xy)|=|\tr\rho_i(xy)|$, then $\rho_{i'}$ is conjugate to $\rho_i$ or to $\epsilon\rho_i$, by the parametrization of the arc.
In the first case $\tr\rho_{i'}([C_\gamma])=\tr\rho_i([C_\gamma])$ holds for every $\gamma$,
and in the second case $\tr\rho_{i'}([C_\gamma])=(-1)^{\Psi(\gamma)}\tr\rho_i([C_\gamma])$ holds for every $\gamma$,
so that the two coordinates are not independent.
\end{rem}

\section{Branched covers and subgroups of finite index}
\label{sec.branched.compact}

The law for compact groups passes to the preimages under finite branched covers.
For a base obeying the law for every compact Lie group with semisimple identity component,
this is \cite[Theorem B]{UekiSC} for the norm ordering, and it is proved in \cite{UekiSC} for every compatible ordering;
the proof induces the representation from the subgroup of finite index that corresponds to the cover
and applies the law for the base to the closure of the image of the induced representation.
We present the argument for the modular knots and an arbitrary compact Lie group,
in the notation of \Cref{subsec.setting.compatible};
the law for the base is now \Cref{thmA.compact} \textup{(1)},
which applies to every compact Lie group, so the semisimplicity is not needed anywhere.

\begin{thm} \label{thm.branched.compact}
Let $\varpi\colon M_\varpi\to S^{3}$ be a finite cover branched along a finite sublink of $\mca{K}$,
the union of the modular knots and of $K_{2,3}$.
Then every ordering of the components of $\varpi^{-1}(\mca{K})$ that is compatible with the trace ordering of $\mca{K}$ obeys the Chebotarev law for every compact Lie group in $M_\varpi$,
the law being defined by \eqref{eq.law.compact} with $\pi_1(M_\varpi-L_n)$ in place of $\pi_1(S^{3}-L_n)$.
The norm ordering is compatible, and so is the ordering by the periods of the lifted modular flow,
with the components of $\varpi^{-1}(K_{2,3})$ placed at fixed positions.
\end{thm}

\begin{proof}
Let $(K'_\nu)_\nu$ be a compatible ordering, let $n\in\Z_{\geq 1}$,
put $L'=\bigcup_{\nu\leq n}K'_\nu$, let $G$ be a compact Lie group,
let $\rho\colon\pi_1(M_\varpi-L')\to G$ have dense image, and let $F$ be a continuous class function on $G$.
By the Peter--Weyl theorem we may assume that $F$ is the character $\chi_V$ of an irreducible unitary representation $V$ of $G$.
Choose $m$ such that $L_m=\bigcup_{i\leq m}K_i$ contains the branch locus and $\varpi(L')$,
put $\Pi=\pi_1(S^{3}-L_m)$ and $H=\pi_1(M_\varpi-\varpi^{-1}(L_m))$,
a subgroup of finite index of $\Pi$ since $\varpi$ restricts to a cover,
and write $\rho$ also for the composite of $H\to\pi_1(M_\varpi-L')$, which is surjective, with $\rho$.
Let $\rho_V=V\circ\rho$, let $\operatorname{Ind}\rho_V$ denote the representation of $\Pi$ induced from $\rho_V$,
of dimension $\deg\varpi\cdot\dim V$, and let $G'_V$ denote the closure of its image.

Let $j>m$, and let $\gamma\in\Pi$ represent the class of $K_j$.
The components of $\varpi^{-1}(K_j)$ correspond to the orbits of the right multiplication by $\gamma$ on $H\backslash\Pi$,
a component of degree $k$ corresponding to an orbit with $k$ elements,
and for a fixed point $Hg$ the element $g\gamma g^{-1}$ lies in $H$ and represents the class of the corresponding component in $\pi_1(M_\varpi-L')$.
Hence the sum of $\chi_V\circ\rho$ over the components of degree one of $\varpi^{-1}(K_j)$ is
\[\sum_{Hg\gamma=Hg}\chi_V\bigl(\rho(g\gamma g^{-1})\bigr)
=\chi_{\operatorname{Ind}\rho_V}(\gamma),\]
by the formula for the character of an induced representation.
We apply \Cref{thmA.compact} \textup{(1)} to the homomorphism $\operatorname{Ind}\rho_V\colon\Pi\to G'_V$,
whose image is dense by the choice of $G'_V$,
and to the class function $\tr$ of $G'_V$, the character of its tautological representation:
the average of $\chi_{\operatorname{Ind}\rho_V}(\gamma_j)$,
where $\gamma_j\in\Pi$ represents the class of $K_j$,
over $m<j\leq J$ converges, as $J\to\infty$, to $\int_{G'_V}\tr(g)\,dg$,
which is the dimension of the invariants of $G'_V$, that is, of $\Pi$, in $\operatorname{Ind}\rho_V$,
and that dimension equals $\dim V^{H}$ by Frobenius reciprocity.
If $V$ is trivial, then the limit is $1$, so that the number $D(J)$ of the components of degree one over $K_{m+1},\dots,K_J$ satisfies $D(J)\sim J$;
this step uses only the law for the finite image of the permutation representation of $\Pi$ on $H\backslash\Pi$, and does not depend on $\rho$.
If $V$ is nontrivial, then $V^{H}=0$, because $\rho(H)$ is dense in $G$ and $V$ is irreducible,
so that the sum of $\chi_V\circ\rho$ over the components of degree one over $K_{m+1},\dots,K_J$ is $o(J)=o(D(J))$.
Hence the average of $F\circ\rho$ over the components of degree one over $K_{m+1},\dots,K_J$ converges to $\int_G F\,dg$.

The number $J(\nu)$ is finite for every $\nu$, because $D(J)\to\infty$ holds as $J\to\infty$.
By \textup{(C1)}, the components of degree one among $K'_{n+1},\dots,K'_\nu$ are,
up to $o(\nu)$ of them, the components of degree one over $K_{m+1},\dots,K_{J(\nu)}$,
the finitely many components over $K_1,\dots,K_m$ being neglected;
by \textup{(C2)}, the components of higher degree among $K'_{n+1},\dots,K'_\nu$ are $o(\nu)$ in number.
Since $J(\nu)\to\infty$ and $D(J)\sim J$, this implies $J(\nu)\sim\nu$,
and the average of $F\circ\rho$ over $K'_{n+1},\dots,K'_\nu$ converges to $\int_G F\,dg$, since $F$ is bounded.

It remains to verify that the two orderings are compatible.
Let $\deg\varpi=d$.
For the norm ordering, let $X$ be the norm of $K'_\nu$.
Every component of degree one over $K_1,\dots,K_{X-1}$ has norm at most $X-1$,
hence is among $K'_1,\dots,K'_\nu$, so $J(\nu)\geq X-1$,
and the components of degree one among $K'_1,\dots,K'_\nu$ over some $K_i$ with $i>J(\nu)$ have $i=X$, so they are at most $d$ in number;
a component of degree $k\geq 2$ with norm at most $X$ lies over some $K_i$ with $i\leq X^{1/2}$,
so those components are at most $dX^{1/2}$ in number, while $\nu\geq D(X-1)\sim X$.
For the ordering by the periods, let $T$ denote the period of $K'_\nu$,
so that $K'_1,\dots,K'_\nu$ have periods at most $T$,
and let $\pi(T)\sim e^{T}/T$ denote the number of the modular knots of length at most $T$ \PI{thm.arith}.
Every component of degree one over a knot of length less than $T$ has period less than $T$, hence is among $K'_1,\dots,K'_\nu$;
since the trace ordering is by length, $J(\nu)$ is at least the number of the knots of length less than $T$,
and the components of degree one among $K'_1,\dots,K'_\nu$ over some $K_i$ with $i>J(\nu)$ lie over knots of length exactly $T$,
which are $o(\pi(T))$ in number by the argument of \PI{lem.ties},
as in the proof of \PI{thm.arith}, while $\nu\geq D(J(\nu))\sim\pi(T)$;
the finitely many components of $\varpi^{-1}(K_{2,3})$ change $J(\nu)$ by a bounded amount.
A component of degree $k\geq 2$ with period at most $T$ lies over a knot of length at most $T/2$,
so those components are at most $d\,\pi(T/2)$ in number, which is $o(\pi(T))$.
\end{proof}

\begin{cor} \label{cor.finiteindex.compact}
Let $\Gamma'\subset \PSL_2\Z$ be a subgroup of finite index, let $M'=T^{1}(\Gamma'\backslash\bb{H}^{2})$,
and let $(K_i)_i$ denote the closed orbits of the geodesic flow on $M'$, ordered by hyperbolic length.
Then $(K_i)_i$ obeys the Chebotarev law for every compact Lie group in $M'$.
In particular, for every homomorphism $\rho\colon \pi_1(M')\to \mathrm{SU}(2)$ with dense image,
the values $\frac12\tr\rho([K_i])$ are equidistributed with respect to $\frac{2}{\pi}\sqrt{1-x^{2}}\,dx$.
\end{cor}

\begin{proof}
As in the proof of \PI{thm.finiteindex}, $M'$ is the complement of the preimage of $K_{2,3}$ in a finite cover $M_\varpi\to S^{3}$ branched along $K_{2,3}$,
the closed orbits of the geodesic flow on $M'$ are the components of the preimages of the modular knots,
and their hyperbolic lengths are their periods.
\Cref{thm.branched.compact} for the ordering by the periods,
with the components of the preimage of $K_{2,3}$ placed first,
which changes the ordering by finitely many positions, is the assertion,
since $\pi_1(M'-L)=\pi_1(M_\varpi-(L\cup\varpi^{-1}(K_{2,3})))$ for every finite union $L$ of closed orbits.
\end{proof}

For a generic ordering of the base, the law for the cover needs the law for the base with respect to the induced representations,
and \Cref{thmA.general} \textup{(1)} makes that condition explicit.
Let $\varpi\colon M_\varpi\to S^{3}$ be a finite cover branched along a finite union of modular knots and of $K_{2,3}$,
let $L'$ be a finite union of components of the preimages, let $G$ be a compact Lie group,
and let $\rho\colon \pi_1(M_\varpi-L')\to G$ have dense image.
Let $\Pi$ denote the fundamental group of the complement in $S^{3}$ of a finite union $L_m$ of modular knots
and $K_{2,3}$ that contains the branch locus and $\varpi(L')$,
let $H\subset \Pi$ denote the subgroup of finite index corresponding to the cover,
and for an irreducible unitary representation $V$ of $G$ let $G'_V$ denote the closure of the image of the induced representation $\operatorname{Ind}_H^{\Pi}(V\circ\rho)$.

\begin{prop} \label{prop.branched.generic}
Let $h$ be a non-arithmetic roof function, and assume that, for every nontrivial irreducible $V$,
no pair of a nontrivial continuous character $\chi'$ of $G'_V$
and $a\in \R$ satisfies $\chi'\bigl(\operatorname{Ind}_H^{\Pi}(V\circ\rho)([C_\gamma])\bigr) =e^{ia\ell_h(\gamma)}$ for all but finitely many modular knots $C_\gamma\not\subset L_m$.
Then every ordering of the components of the preimages that is compatible with the ordering of the modular knots by $\ell_h$ obeys the equality \eqref{eq.law.compact} for $\rho$ and every continuous class function of $G$.
The assumption holds if $G^{0}$ is semisimple.
For the ordering by the hyperbolic length in place of $\ell_h$,
the conclusion holds for every $G$ without assumption, by \Cref{thm.branched.compact}.
\end{prop}

\begin{proof}
By \Cref{thmA.compact} \textup{(2)} and the paragraph after \Cref{thmA.compact},
the assumption says that the modular knots obey the law for $(G'_V,\operatorname{Ind}_H^{\Pi}(V\circ\rho))$ for every nontrivial irreducible $V$,
which is what the proof of \Cref{thm.branched.compact} uses for the base;
for the trivial $V$, which that proof uses to count the components of degree one,
the group $G'_V$ is the finite image of the permutation representation of $\Pi$ on $H\backslash\Pi$,
and the law holds for it by \PI{thm.modular}.
If $G^{0}$ is semisimple, then so is $(G'_V)^{0}$, by \cite{UekiSC},
and a character of finite image cannot satisfy the displayed equality with $a\neq 0$,
by \PI{lem.refined.nonarithmetic},
while $a=0$ is excluded by \Cref{lem.generate.linking} for $L=L_m$ and by the density of the image, as in the proof of \Cref{lem.cusp}.
\end{proof}

\begin{rem} \label{rem.branched.generic}
The assumption of \Cref{prop.branched.generic} is stated for the induced representations rather than for $(G,\rho)$ itself,
because a character of $G'_V$ is not necessarily the pullback of a character of $G$:
the determinant of the induced representation is the product of the $\dim V$-th power of the sign of the permutation representation of $\Pi$ on $H\backslash\Pi$
and of $\det(V\circ\rho)$ composed with the transfer $\Pi\to H^{\rm ab}$,
and the closure $G'_V$ may have further characters.
A description of the exceptions for the cover in terms of $(G,\rho)$ and of the transfer is left open.
\end{rem}

\subsection*{Acknowledgments}
The author has been partially supported by JSPS KAKENHI Grant Numbers JP19K14538 and JP23K12969.

This work was carried out between August and November 2026 in extensive collaboration with Claude (Anthropic; models Fable 5, 5.1, and 5.2 and Opus 5.5, accessed through claude.ai) on the mathematics, on the numerical experiments, and on the manuscript.
All mathematical content has been verified by the author, who takes full responsibility for it.

{\small 
\bibliographystyle{amsalpha}
\bibliography{ju_mkc2_Lie_arXiv.bbl
}

\medskip
\noindent Department of Mathematics, Ochanomizu University,
2-1-1 Otsuka, Bunkyo-ku, Tokyo 112-8610, Japan

\noindent\textit{E-mail address}: \texttt{uekijun46@gmail.com}
\end{document}